\documentclass[12pt]{amsart}
\usepackage{amssymb,amsmath, amsthm,latexsym}
\usepackage{a4wide}

\theoremstyle{plain}

\newtheorem{lem}{Lemma}[section]
\newtheorem{theo}[lem]{Theorem}
\newtheorem{prop}[lem]{Proposition}
\newtheorem{corollary}[lem]{Corollary}

\parskip=\medskipamount
\font\k=cmr7

  \newcommand {\inte}{\mbox{\k int}}
  \newcommand {\ext}{\mbox{\k ext}}
  \newcommand {\loc}{\mbox{\k loc}}
  \newcommand {\C}{{\mathbb C}}
  \newcommand {\N}{{\mathbb N}}
  \newcommand {\R}{{\mathbb R}}
  \newcommand {\Z}{{\mathbb Z}}

\renewcommand {\H}{{\mathcal H}}

  \newcommand {\E}{{\mathcal E}}

 \newcommand {\ov}{\overline}

\renewcommand{\Re}{\operatorname{Re}}
\newcommand{\Tr}{\operatorname{Tr}}

\newcommand{\Id}{\operatorname{Id}}
\newcommand{\Hom}{\operatorname{Hom}}

\newcommand{\GL}{\operatorname{GL}}

\begin{document}
\title[Decomposition of determinants]
{\bf Regularized
 determinants of Laplace type operators, analytic surgery and relative
 determinants} 
\date{\today}

\author{J\"orn M\"uller}
\author{Werner M\"uller}
\address{Universit\"at Bonn\\
Mathematisches Institut\\
Beringstrasse 1\\
D -- 53115 Bonn, Germany}
\address{Universit\"at Bonn\\
Mathematisches Institut\\
Beringstrasse 1\\
D -- 53115 Bonn, Germany}
\email{joern@uni-bonn.de}
\email{mueller@math.uni-bonn.de}


\begin{abstract}Let $M$ be a compact Riemannian manifold in which $Y$ is an
embedded hypersurface separating $M$ into two parts. Assume that the metric is
a product on a tubular neighborhood $N$ of $Y$. Let $\Delta$ be a Laplace type
operator on $M$ adapted to the product structure on $N$. Under certain
 additional assumptions on $\Delta$, we
 establish an
asymptotic expansion for the logarithm of the
 regularized determinant $\det\Delta$ of $\Delta$
if   the tubular neighborhood $N$ is stretched to a cylinder of infinite
length. We use the asymptotic expansions to derive adiabatic splitting formulas
for regularized determinants.
\end{abstract}
\maketitle

\section{Introduction}

In this paper we study the behaviour of regularized determinants of Laplace
type operators 
with respect to  certain  singular deformations which are related to
 {\it analytic surgery}.
Analytic surgery is a method developed by Mazzeo and Melrose \cite{MM} to
study the behaviour of global spectral invariants of Dirac- and Laplace
operators with respect to decompositions of the underlying Riemannian
manifolds.  

The singular deformations  that we consider
in this paper are defined in the following way. Let $M$ be a compact Riemannian
manifold and let $Y$ be an embedded hypersurface
in $M$ such that $M-Y$ consists of two components $M_1$ and $M_2$. Assume that
the metric in a collar neighborhood $N$ of $Y$ is a product. Then by ''analytic
surgery'' we mean the stretching of the collar neighborhood $N$ to a cylinder
of infinite length. In this way we get a family of Riemannian manifolds 
$(M_r,g_r)$, $r\ge 1$. The singular limit of this family is the disjoint union
of two 
manifolds with cylindrical ends $M_{1,\infty}$ and $M_{2,\infty}$. Let
$\Delta\colon C^\infty(M,E)\to C^\infty(M,E)$ be a Laplace type operator on
$M$ which is adapted to the product structure on $N$. Then we can define an
associated family of Laplace type operators $\Delta_r$ on $M_r$ and the main
purpose of this paper is to study the  behaviour of $\det(\Delta_r)$
as $r\to\infty$. Under some additional assumptions on $\Delta$ we will show
that $\log\det\Delta_r$ has an asymptotic 
expansion and the main ingredient of the constant term of this expansion are
the relative determinants $\det(\Delta_{i,\infty},\Delta_0)$, $i=1,2$,
associated to the Laplacian on the manifolds with cylindrical ends
$M_{1,\infty}$ and $M_{2,\infty}$, respectively. Here the relative
determinants are defined as in \cite{Mu1}.  For surfaces our results are
related to the work of Bismut and Bost  \cite{BB} who studied the Quillen
metric on the determinant line bundle associated to a family of complex
curves with singular fibers. 

We also consider the analogous problem  for a compact
manifold with boundary where we stretch a collar neighborhood of the
boundary to an  infinite half-cylinder. The singular limit of the associated
family  of Riemannian manifolds with boundary is a manifold  with a
cylindrical end. The relative determinant of the Laplacian on the manifold
with cylindrical end arises in the same manner as above in the asymptotic 
expansion of the determinant of the Dirichlet Laplacians. 
This gives a new interpretation of relative determinants. 

If we compare the asymptotic expansions of the 
determinants of the Laplacians on the manifolds $M_r$, $M_{1,r}$ and
$M_{2,r}$, respectively, 
obtained by  stretching the corresponding collar neighborhoods of $Y$, we
recover the adiabatic decomposition formulas of Park and Wojciechowski 
\cite{PW1}, \cite{PW3}. We also establish a gluing formula for relative
determinants of Laplace type operators on manifolds with cylindrical ends.

In the present paper we consider Laplace operators of two types. First we
assume that the induced Laplace operator on $Y$ is invertible and that
the Laplacians $\Delta_{i,\infty}$, $i=1,2$,  on $M_{i,\infty}$ have no nonzero
$L^2$-solutions. This simplifies the constructions. The second case that we
consider are Bochner-Laplace operators. In a followup paper we will study the
case of Dirac-Laplace operators $\Delta=D^2$.

Now we describe the content of the paper in more detail. Let $(X,g)$ be a
Riemannian manifold and let $E\to X$ be a Hermitian vector bundle. First
recall that a Laplace type operator
$$\Delta\colon C^\infty(X,E)\to C^\infty(X,E)$$
is a second order elliptic differential operator which is symmetric,
nonnegative and whose principal symbol is given by 
$$\sigma_{\Delta}(x,\xi)=\parallel\xi\parallel^2\Id_{E_x}.$$
Suppose that $X$ is a compact manifold with boundary $\partial X$, which may
be empty. We impose Dirichlet boundary conditions at $\partial X$ and denote
the corresponding selfadjoint extension by $\Delta_{D}$. This is a
selfadjoint nonnegative operator in $L^2(X,E)$. The regularized determinant
$\det\Delta_{D}$ of $\Delta_{D}$ is defined in the usual way by
$$\det\Delta_{D}=\exp\left(-\frac{d}{ds}
\zeta_{\Delta_{D}}(s)\Big|_{s=0}\right),$$ 
where $\zeta_{\Delta_{D}}(s)$ is the zeta function of $\Delta_{D}$.

Our first result is a gluing formula for relative determinants of Laplace
type operators on a manifold $X$ with a cylindrical end.
 By definition, $X$ has a decomposition 
$$X=M\cup_Y Z,\quad Z=\R^+\times Y,$$
where $M$ is a compact manifold with boundary $Y$ and the metric $g^X$ of $X$ 
is a product on $\R^+\times Y$. Let $E\to X$ be a hermitian vector bundle.  
 We assume that there exist a hermitian vector bundle $E_0\to Y$ such that 
$E|Z\cong \mathrm{pr}_Y^*E_0$ and that the fiber metric $h^E$ of $E$ is a
product  
on $\R^+\times Y$.  Let $\Delta\colon C^\infty(X,E)\to C^\infty(X,E)$
be a Laplace type operator on $X$.  We assume that the restriction of
$\Delta$ to $Z$ satisfies 
\begin{equation}\label{1.1a}
\Delta|_Z=-\frac{\partial^2}{\partial u^2}+\Delta_Y,
\end{equation}
where $\Delta_Y$ is a Laplace type operator on $Y$. This implies that
 $\Delta_X$ is essentially selfadjoint in 
$L^2$. We will  denote the unique selfadjoint extension of $\Delta_X$ by the
same letter. Consider the operator
$$-\frac{\partial^2}{\partial u^2}+\Delta_Y\colon C^\infty_c(\R^+\times Y,E)
\to L^2(\R^+\times Y,E).$$
and impose Dirichlet boundary conditions at $\{0\}\times Y$. Let
$\Delta_0$ be the corresponding selfadjoint extension.  Then
$\Delta,\Delta_0)$
is a pair of self-adjoint operators which satisfies  conditions 1)--3) in
\cite[p.312]{Mu1} which are needed to define 
the relative regularized determinant $\det(\Delta,\Delta_0)$. 

Let $\Delta_M$ denote the restriction of $\Delta$ to $M$ and let
$\Delta_{M,D}$ be the selfadjoint extension obtained by imposing Dirichlet
boundary conditions at $\partial M$.   We assume that $\Delta_{M,D}$ is 
invertible. This assumption is satisfied in many cases. Suppose, for example, 
 that $D\colon C^\infty(X,E)\to C^\infty(X,E)$ is a 
Dirac  operator and  $\Delta=D^2$. Then it follows from \cite{Ba}
that $\Delta_D$ is invertible. In particular, if 
$\Delta_p\colon \Lambda^p(X)\to\Lambda^p(X)$ is the Laplacian on $p$-forms
on a compact manifold with boundary, then $\Delta_{p,D}$ is invertible.
Other examples are Bochner-Laplace operators.

If $\Delta_{M,D}$ is invertible, then the Dirichlet-to-Neumann operator $R$ 
with respect to the hypersurface $Y\cong \{0\}\times Y\subset X$ can be
defined in the usual way.
This is a pseudo-differential operator of order 1 on $Y$  which 
is selfadjoint and nonnegative. So $R$  has a well-defined 
determinant $\det R$. 

The last ingredient of the gluing formula is defined in terms of the space
$\H$ of extended $L^2$-solutions of $\Delta$. Recall that a section
$\varphi\in C^\infty(X,E)$ is called an 
extended  $L^2$-solution of $\Delta$, if $\varphi$ is a bounded solution of
$\Delta\varphi=0$
and its restriction to  $\R^+\times Y$  has the form 
$$\varphi(u,y)=\phi(y)+\psi(u,y),$$
where $\psi$ is in $L^2$ and $\phi\in\ker\Delta_Y$. In this case, $\phi$ is
called the limiting value of $\varphi$. Let $V^+\subset \ker\Delta_Y$ be the
space  of all limiting values of extended $L^2$-solutions of $\Delta$.
Given $\phi\in V^+$, let $E(\phi,\lambda)$ be the
associated generalized eigensection of $\Delta$ (cf. \cite{Mu4}). Then  
$E(\phi,\lambda)$ is holomorphic at $\lambda=0$ and $E(\phi,0)$ is an
extended $L^2$ solution of $\Delta$ with limiting value $2\phi$. Let
$\rho_Y\colon C^\infty(X,E) \to C^\infty(Y,E|Y)$ denote the
restriction map and set $\H_Y:=\rho_Y(\H)$. We show that $\rho_Y:\H\to \H_Y$
is an 
isomorphism.  Let $\varphi_1,...,\varphi_k$
be an orthonormal basis of $\ker \Delta$ and let $\phi_1,...,\phi_l$ be
an orthonormal basis of $V^+$. Put $\psi_i=\rho_Y(\varphi_i)$, if
$1\le i\le k$, and $\psi_{k+j}=\frac{1}{2}\rho_Y(E(\phi_{j},0))$, if $1\le j\le
l$. Put  $a_{ij}=\langle\psi_i,\psi_j\rangle_Y$, $1\le i,j\le k+l$ and let 
$A$ be the $(k+l)\times(k+l)$-matrix with entries $a_{ij}$. We are now ready
to state our first main result which is the following theorem.

\begin{theo}\label{th1.1a}
Assume that $\Delta_{M,D}$ is invertible. 
Let $h_Y=\dim\ker \Delta_Y$ and denote by $\zeta_Y(s)$ is the zeta function of 
$\Delta_Y$. Then
$$\frac{\det(\Delta,\Delta_0)}{\det(\Delta_{M,D})}=2^{-\zeta_Y(0)-h_Y}
\frac{\det R}{\det A}.$$
\end{theo}
The same result has been proved independently by Loya and Park \cite{LP}.

Now assume that $(M,g)$ is an oriented closed connected $n$-dimensional
Riemannian  
manifold and let $Y$ be a hypersurface of $M$ such that $M-Y$ consists of two
components. We denote the closure of the components of $M-Y$ by $M_1$ and
$M_2$. Thus $M_1$ and $M_2$ are compact manifolds with common boundary $Y$
such that
\begin{equation}\label{1.1c}
M=M_1\cup_Y M_2,\quad \partial M_1=\partial M_2=Y.
\end{equation}
Let $E\to M$ be a Hermitian vector bundle and let
$$\Delta_M\colon C^\infty(M,E)\to C^\infty(M,E)$$
be a Laplace type operator.
We assume that there exists a tubular neighborhood $N$ of $Y$ which is
diffeomorphic to $[-1,1]\times Y$ such that all geometric structures are
products over $N$, i.e., $g|_N=du^2+g^Y$, there exists a
Hermitian vector bundle $E_0\to Y$ such that $E|N=\mathrm{pr}^*_Y(E_0)$ and 
\begin{equation}\label{1.1}
\Delta_M\big|_N=-\frac{\partial^2}{\partial u^2}+\Delta_Y,
\end{equation}
where $\Delta_Y\colon C^\infty(Y,E_0)\to C^\infty(Y,E_0)$ is a Laplace type
operator on $Y$. Let $\Delta_{M_i}$ be the restriction of $\Delta_M$ to
$M_i$, $i=1,2$. We assume that $\Delta_{M_1,D}$ and $\Delta_{M_2,D}$ are
invertible (see the above remark).

We define a family of Riemannian manifolds $(M_r,g_r)$, $r>0$, as follows.
Given $r>0$, let $N_r=[-r,r]\times Y$ and set
\begin{equation}\label{1.4b}
M_r=M_1\cup_Y N_r\cup_Y M_2,
\end{equation}
where $\partial M_1$ is identified with $\{-r\}\times Y$ and $\partial M_2$ 
with $\{r\}\times Y$. Since $g$ is a product in a neighborhood of $Y$, it has 
a canonical extension to a metric $g_r$ on $M_r$ such that
$g_r|_{N_r}=du^2+g^Y$.  Similarly, $E\to M$ and 
$\Delta_M$ have natural extensions $E_r\to M_r$ and $\Delta_{M_r}$ to $M_r$.
Our main purpose is to study the asymptotic behavior of $\det(\Delta_{M_r})$
as $r\to\infty$.  To describe the result
we need some more notation. Set
$$M_{i,\infty}=M_i\cup_Y (\R^+\times Y),\;\;i=1,2.$$
This is a manifold with a cylindrical end $Z=\R^+\times Y$. The disjoint 
union of $M_{1,\infty}$
and $M_{2,\infty}$ may be regarded as the singular limit of $M_r$ as 
$r\to\infty$. Let $\Delta_{i,\infty}$ be the canonical extension of 
$\Delta_M\big|_{M_i}$ to $M_{i,\infty}$ which is defined by
$$\Delta_{i,\infty}\big|_{M_i}=\Delta_M\big|_{M_i},\quad
\Delta_{i,\infty}\big|_{\R^+\times Y}=-\frac{\partial^2}{\partial u^2}+
\Delta_Y.$$
Then $\Delta_{i,\infty}$ is essentially selfadjoint in $L^2$. We denote the
unique selfadjoint extension of $\Delta_{i,\infty}$ by the same letter. Let
$\Delta_0$ be as in Theorem \ref{th1.1a} and let 
$\det(\Delta_{i,\infty},\Delta_0)$ be the relative determinant  \cite{Mu1}.

 Let
\begin{equation}\label{1.4a}
\xi_Y(s):=\frac{\Gamma(s-1/2)}{\sqrt{\pi}\Gamma(s)}\zeta_Y(s-1/2),
\end{equation}
where $\zeta_Y(s)$ is the zeta function of $\Delta_Y$.
Our first  result  concerning the asymptotic behaviour of the determinant of
a Laplace type operator is obtained under the assumption that all involved
operators are invertible.

\begin{theo}\label{th1.1}
Suppose that $\ker\Delta_Y=\{0\}$ and $\ker\Delta_{i,\infty}=\{0\}$, $i=1,2$.
Then $\Delta_{M_r}$ is invertible for $r\ge r_0$ and
\begin{equation}\label{1.2}
\lim_{r\to\infty}e^{r\xi'_Y(0)}\det\Delta_{M_r}
= (\det\Delta_Y)^{-1/2}
\prod_{i=1}^2\det(\Delta_{i,\infty},\Delta_0).
\end{equation}
\end{theo}
In particular, the assumption of Theorem \ref{th1.1} are satisfied for the
operator
$\Delta_M+\lambda$, where $\lambda>0$. Let $\xi_Y(s,\lambda)$ be defined as in
(\ref{1.4a}) with $\zeta_Y(s)$ replaced by the zeta function
$\zeta_Y(s,\lambda)$ of $\Delta_Y+\lambda$. Then we get
\begin{corollary}\label{c1.3a}
Let $\lambda>0$. Then
\begin{equation*}
\begin{split}
\lim_{r\to\infty} e^{r\xi_Y^\prime(0,\lambda)}\det(\Delta_{M_r}+\lambda)
=\det(\Delta_Y+\lambda)^{-1/2}
\prod_{i=1}^2\det(\Delta_{i,\infty}+\lambda,\Delta_0+\lambda).
\end{split}
\end{equation*}
\end{corollary}

We note that (\ref{1.2})  also holds if $M$ has a nonempty
boundary $\partial M$. In this case we impose Dirichlet boundary conditions at
$\partial M$. 

In particular, we may consider a separating hypersurface which is 
parallel to the boundary. This is a special case which we consider separately.
Let $X_0$ be a compact manifold with boundary $Y$ and assume that all geometric
structures are products in a collar neighborhood of $Y$. Let
$X_r=X_0\cup_Y([0,r]\times Y)$ and let $\Delta_{X_r,D}$ be the selfadjoint
extension of the corresponding Laplace operator with respect to Dirichlet
boundary conditions. Then the analogous statement to Theorem \ref{th1.1} is
\begin{prop}\label{p1.4}
Assume that  $\Delta_Y$ and $\Delta_\infty$ 
are invertible. Then $\Delta_{X_r,D}$ is invertible for $r\ge r_0$ and
\begin{equation}\label{1.7}
\lim_{r\to\infty}e^{r\xi_Y^\prime(0)/2}\det\Delta_{X_r,D}=(\det\Delta_Y)^{-1/2}
\det(\Delta_{\infty},\Delta_0),\quad r\to\infty.
\end{equation}
\end{prop}
{\bf Remark.} Under the same assumptions as in Theorem \ref{th1.1} and
Proposition \ref{p1.4}, respectively,  Lee \cite{Le4} has also obtained 
asymptotic expansions for $\log\det\Delta_{M_r}$ and $\log\det\Delta_{X_r,D}$,
which are different from ours. The relation between \cite{Le4} and our results
is given by Theorem \ref{1.1a}.

Especially consider the manifolds with boundary $M_1$ and $M_2$ of the
decomposition (\ref{1.1c}) of $M$. 
Let $M_{i,r}=M_i\cup_Y([0,r]\times Y)$, $i=1,2$. If we apply (\ref{1.7}) 
to $\Delta_{M_{i,r},D}$ and compare it to (\ref{1.2}) , 
we obtain
\begin{corollary}\label{c1.5}
Let $M$ be closed. Assume that  $\Delta_Y$ and 
$\Delta_{i,\infty}$, $i=1,2$, are invertible. Then $\Delta_{M_r}$ and
$\Delta_{M_{i,r},D}$, $i=1,2$,  are invertible for $r\ge r_0$ and
\begin{equation}\label{1.8}
\lim_{r\to\infty}\frac{\det\Delta_{M_r}}{\det\Delta_{M_{1,r},D}
\det\Delta_{M_{2,r},D}}=(\det\Delta_Y)^{1/2}.
\end{equation}
\end{corollary}
This is the ''adiabatic decomposition formula'' established by Park and
Wojciechowski in  \cite{PW1}. 

Next we study the case of a Bochner-Laplace operator. 
Let  $\nabla$ is a metric connection on $E$ which is a product on 
$N$. Let $\Delta_M=\nabla^*\nabla$ be the associated Bochner-Laplace operator.
Then $\nabla$ has  canonical extensions to a connection $\nabla^r$ on 
$E_r\to M_r$ and $\nabla^{i,\infty}$ on $E_{i,\infty}\to M_{i,\infty}$,
respectively, and $\Delta_{M_r}$  and $\Delta_{i,\infty}$ are the 
corresponding Bochner-Laplace operators. We need to introduce some further
notation. Let 
$$S_i(0)\colon\ker\Delta_Y\to\ker\Delta_Y,\quad i=1,2,$$
 denote the on-shell scattering operator at energy zero
 associated to $(\Delta_{i,\infty},\Delta_0)$ (see e.g. \cite{Mu4}). 
This operator satisfies $S_i(0)^2=\Id$. Let 
$$\ker\Delta_Y=V_i^+\oplus V_i^-,\quad i=1,2,$$
be the decomposition of $\ker\Delta_Y$ into the $\pm1$-eigenspaces of
 $S_i(0)$. Let $C_{12}$ denote the restriction of $S_1(0)S_2(0)$ to the
 orthogonal complement of 
$(V_1^+\cap V_2^+)\oplus(V_1^-\cap V_2^-)$ in $\ker\Delta_Y$.
 Then our next result is  

\begin{theo}\label{th1.6}
Let $\Delta_M=\nabla^*\nabla$ be a Bochner-Laplace operator. Let $h_Y=
\dim\ker\Delta_Y$ and  $h=\dim V_1^++\dim V_2^+ -2\dim V_1^+\cap V_2^+$.  Then 
\begin{equation}\label{1.9}
\begin{split}
\lim_{r\to\infty}r^{h-h_Y}e^{r\xi^\prime_Y(0)}\det\Delta_{M_r}=&
2^{2h_Y-h}(\det\Delta_Y)^{-1/2}\\
&\cdot\det\left((\Id-C_{12})/2\right)
\prod_{i=1}^2\det(\Delta_{i,\infty},\Delta_0).
\end{split}
\end{equation}
\end{theo}

If we specialize Theorem \ref{th1.6} to the case of the Laplacian 
$\Delta=d^*d$  on functions on a closed surface $M$, 
we obtain 
$$\det\Delta_r\sim 2\det(\Delta_{1,\infty},\Delta_0)
\det(\Delta_{2,\infty},\Delta_0)re^{-\pi r/3}$$
as $r\to\infty$. This is Theorem 13.7 of \cite{BB} with an explicit constant 
expressed in terms of relative determinants.

As in Proposition \ref{p1.4}, we may also consider the case of a compact 
Riemannian manifold $X_0$ with boundary $Y$. For a Bochner-Laplace operator
 on $X_0$ it follows from \cite{Ba} that $\Delta_{X_r,D}$ is invertible.
Let $V^+\subset \ker\Delta_Y$ be the $+1$-eigenspace of the scattering 
operator $S(0)$ and let $h^+=\dim V^+$. The analogous result to (\ref{1.7})
is 
\begin{equation}\label{1.10}
\lim_{r\to\infty}r^{h^+ -h_Y}e^{r\xi_Y^\prime(0)/2}\det\Delta_{X_r,D}=
2^{h_Y}\left(\det\Delta_Y\right)^{-1/2}\det(\Delta_\infty,\Delta_0).
\end{equation}
Now we apply this again to the manifolds with boundary $M_1$ and $M_2$ of
the decomposition (\ref{1.1c}) of $M$ and compare it to (\ref{1.9}). In
this way we get
\begin{theo}\label{th1.7}
Let the notation be as in Theorem \ref{th1.6} and let $h_{12}=\dim V_1^+\cap
V_2^+$.Then 
$$\lim_{r\to\infty}\frac{r^{h_Y-2h_{12}}\det\Delta_{M_r}}
{\det\Delta_{M_{1,r},D}\det\Delta_{M_{2,r},D}}=2^{-h}
\left(\det\Delta_Y\right)^{1/2}\det\left((\Id-C_{12})/2\right).$$
\end{theo}

{\bf Remark.} This result was first proved by by Park and Wojciechowski 
\cite{PW3} under an additional assumption, called Condition A \cite[p.4]{PW3},
which rules out the existence of exponentially decreasing eigenvalues of
$\Delta_{M_r}$. As pointed out by Park and Wojciechowski, their assumption
implies that 1 is not an eigenvalue of $S_1(0)S_2(0)$. This has the
consequence that $V_1^\pm\cap V_2^\pm=\{0\}$, which in turn implies that
  $h=h_Y$ and $h_{12}=0$ and Theorem
\ref{th1.7} specializes to Theorem 0.1 of \cite{PW3}.

Next we explain some of the main ideas of the proofs.
The strategy to prove Theorem \ref{th1.1a} is analogous to the proof of the 
surgery formula in \cite{HZ}. Let $z\in\C-\R_-$. Then the relative
determinant $\det(\Delta+z,\Delta_0+z)$ and the determinant 
$\det(\Delta_{M,D}+z)$ are defined. Moreover the Dirichlet-to-Neumann operator
$R(z)$ with respect to $\Delta+z$ and the hypersurface $Y\subset X$ exists
and the determinant $\det R(z)$ can be defined. Then by Theorem 4.2 of 
\cite{Ca} there is a polynomial $P(z)$ with real coefficients of degree $\le
(n-1)/2$ such that
\begin{equation}\label{1.6}
\frac{\det(\Delta+z,\Delta_0+z)}{\det(\Delta_{M,D}+z)}
=e^{P(z)}\det (R(z)).
\end{equation}
Both sides of this equality have an expansion in $z$ as $z\to0$. We determine 
these expansions and compare the constant terms. This proves Theorem 
\ref{th1.1a}.

To prove Theorems \ref{th1.1} and  \ref{th1.6}, we apply the 
Mayer-Vietoris formula of \cite{BFK}  to $\det(\Delta_{M_r}+\lambda)$, 
$\lambda>0$, with respect to the decomposition (\ref{1.4b})
and take the limit $\lambda\to 0$. To this end we assume that 
$\Delta_{M_1,D}$
and $\Delta_{M_2,D}$ are invertible. Under this assumption the 
Dirichlet-to-Neumann operator $R_r$ with respect to the hypersurface 
$\Sigma_r:=(\{-r\}\times Y)\sqcup(\{r\}\times Y)$ exists and we get a
splitting formula 
for $\det\Delta_{M_r}$. We compare  
this splitting formula with the splitting formulas for 
$\det(\Delta_{i,\infty},\Delta_0)$ given by Theorem \ref{th1.1a}. 
Finally we study the limit of
$\det R_r$ as $r\to\infty$ and compare it to $\det R_{1,\infty}\det
R_{2,\infty}$.  Let $\Delta_{N_r,D}$ be the Laplace operator on $N_r$ with
Dirichlet boundary conditions. Under the assumptions of Theorem \ref{th1.1} or
\ref{th1.6} the limit as $r\to\infty$ of 
$r^{h}\det\Delta_{M_r}(\det\Delta_{N_r,D})^{-1}$ 
exists and 
$$\lim_{r\to\infty}\frac{r^{h}\det\Delta_{M_r}}{\det\Delta_{N_r,D}}
=2^{-h}\det\left((\Id-C_{12})/2\right)
\prod_{i=1}^2\det(\Delta_{i,\infty},\Delta_0).$$
Finally we determine the asymptotic behaviour of $\det\Delta_{N_r,D}$ as
$r\to\infty$. This completes  the proof of Theorem \ref{th1.1} and Theorem
\ref{th1.6}.

\section{Expansion of relative determinants}
\setcounter{equation}{0}
Let $X$ be a manifold with a cylindrical end and let $\Delta$ be a Laplace
type operator on $X$ as above.
In this section we consider the asymptotic expansion of
$\log\det(\Delta+z,\Delta_0+z)$ 
as $z\to0$. We use  the framework introduced in \cite{Mu1}. 
Let $H,H_0$ be two self-adjoint nonnegative linear operators in a separable 
Hilbert space $\H$ such that $e^{-tH}-e^{-tH_0}$ is a trace class operator 
for all $t>0$. Suppose that the following two conditions are satisfied:
\begin{enumerate}
\item[1)] As $t\to0+$, there exists an asymptotic expansion of the form
$$\Tr(e^{-tH}-e^{-tH_0})\sim\sum_{j=0}^\infty a_j t^{\alpha_j},$$
where $-\infty<\alpha_0<\alpha_1<\cdots$ and $\alpha_j\to\infty$. 
\item[2)] There exist $b_0\in\C$, $\rho>0$ such that
$$\Tr(e^{-tH}-e^{-tH_0})\sim b_0+O(t^{-\rho})$$
as $t\to\infty$.
\end{enumerate}
 Set
\begin{equation*}
\begin{aligned}
\zeta_1(s,H,H_0)&=\frac{1}{\Gamma(s)}\int_0^1
t^{s-1}\Tr(e^{-tH}-e^{-tH_0})\;dt,\quad \Re(s)>-\alpha_0;\\
\zeta_2(s,H,H_0)&=\frac{1}{\Gamma(s)}\int_1^\infty
t^{s-1}\Tr(e^{-tH}-e^{-tH_0})\;dt,\quad \Re(s)<0.
\end{aligned}
\end{equation*}
Then $\zeta_1(s,H,H_0)$ admits a meromorphic extension to $\C$ which is
holomorphic at $s=0$. Similarly $\zeta_2(s,H,H_0)$ has a meromorphic 
extension  to the half-plane $\Re(s)<\rho$ which is also holomorphic at $s=0$.
It is given by
\begin{equation}\label{2.1}
\zeta_2(s,H,H_0)=-\frac{b_0}{\Gamma(s+1)}+\frac{1}{\Gamma(s)}
\int_1^\infty t^{s-1}\left(\Tr(e^{-tH}-e^{-tH_0})-b_0\right)\;dt.
\end{equation}
The relative zeta function $\zeta(s,H,H_0)$ is then defined by
$$\zeta(s,H,H_0)=\zeta_1(s,H,H_0)+\zeta_2(s,H,H_0),$$
and the relative determinant by
$$\det(H,H_0)=\exp\left(-\frac{d}{ds}\zeta(s,H,H_0)\big|_{s=0}\right).$$
Let $\lambda> 0$ and define $\det(H+\lambda,H_0+\lambda)$ similarly. 

\begin{prop}\label{p2.1}
As $\lambda\to0+$, we have
$$\log\det(H+\lambda,H_0+\lambda)=b_0\log\lambda+\log\det(H,H_0)+o(1).$$
\end{prop}
\begin{proof}
From the  construction of the analytic continuation of
 $\zeta_1(s,H+\lambda,H_0+\lambda)$ and $\zeta_1(s,H,H_0)$, respectively, 
 it follows immediately that
$$\lim_{\lambda\to0}\zeta_1^\prime(0,H+\lambda,H_0+\lambda)
=\zeta_1^\prime(0,H,H_0).$$
Let $\Re(s)<\rho$. Using (\ref{2.1}) we get 
\begin{equation*}
\begin{aligned}
\zeta_2(s,H+\lambda,H_0+\lambda)&=\frac{1}{\Gamma(s)}\int_1^\infty t^{s-1}
e^{-t\lambda}\Tr(e^{-tH}-e^{-tH_0})\;dt\\
&=\frac{1}{\Gamma(s)}\int_1^\infty t^{s-1}e^{-t\lambda}
\left(\Tr(e^{-tH}-e^{-tH_0})-b_0\right)\;dt\\
&\quad+\frac{b_0}{\Gamma(s)}\int_1^\infty t^{s-1}e^{-t\lambda}\;dt\\
&=\zeta_2(s,H,H_0)+\frac{b_0}{\Gamma(s+1)}+b_0\lambda^{-s}-
\frac{b_0}{\Gamma(s)}\int_0^1 t^{s-1}e^{-t\lambda}\;dt+o(1)
\end{aligned}
\end{equation*}
as $\lambda\to 0+$.
This implies that
$$\zeta_2^\prime(0,H+\lambda,H_0+\lambda)=\zeta_2^\prime(0,H,H_0)
-b_0\log\lambda+o(1)$$
as $\lambda\to0+$.
\end{proof}

In order to apply this result to our case, we need to compute $b_0$. Let 
$\xi(\lambda)$ be the spectral shift function of $(\Delta,\Delta_0)$
\cite[pp. 315]{Mu1}. By (2.16) of \cite{Mu1}, we have
\begin{equation}\label{2.2}
b_0=-\xi(0+).
\end{equation}
So we are reduced to the study of the spectral shift function near zero.
Recall that the spectral shift function is a real valued function in
$L^2_{\loc}(\R)$ which is uniquely determined by the following two properties
\begin{enumerate}
\item $\xi(\lambda)=0$ for all $\lambda<0$.
\item For every $f\in C^\infty_c(\R)$, $f(\Delta)-f(\Delta_0)$ is a trace 
class operator and
$$\Tr(f(\Delta)-f(\Delta_0))=\int_R f^\prime(\lambda)\xi(\lambda)\;d\lambda.$$
\end{enumerate}
Let $\Delta_d$ and $\Delta_{ac}$ denote the restriction of $\Delta$ to the
subspace of $L^2(X,E)$ corresponding to the point spectrum and the absolutely
 continuous spectrum of $\Delta$, respectively. By \cite{Do}, the eigenvalues 
of $\Delta$
have no finite point of accumulation. Hence  
$f(\Delta_d)$ is a trace class operator for every $f\in C^\infty_c(\R)$. 
This implies that $f(\Delta_{ac})-f(\Delta_0)$ is also a trace class operator
for every $f\in C^\infty_c(\R)$. Let $\xi_c(\lambda)$ be the spectral shift
function of $(\Delta_{ac},\Delta_0)$ and let $N(\lambda)$ denote the counting
function of the eigenvalues of $\Delta$. Then it follows from (1) and (2)
that
\begin{equation}\label{2.3}
\xi(\lambda)=-N(\lambda)+\xi_c(\lambda).
\end{equation}
The spectral shift function $\xi_c(\lambda)$ can be determined in the same way
as in Chapter IX of \cite{Mu3}. The manifolds considered in \cite{Mu3}
are manifolds with fibered cusps which are different from the manifolds in
the present paper. However, the structure of the continuous spectrum is
 similar and everything said about the continuous spectrum in 
\cite{Mu3} applies with minor modifications in our case as well.
  Let $\mu_1>0$ be the smallest 
positive eigenvalue of $\Delta_Y$. Let
$$S(s)\colon \ker\Delta_Y\to\ker\Delta_Y,\quad |s|<\sqrt{\mu_1},$$
be the scattering matrix \cite{Mu4}. It is an analytic function. Then it
follows as in the proof of Theorem 9.25 of \cite{Mu3} that
$$\xi_c(\lambda)=-\frac{1}{4}(\Tr(S(0))+\dim\ker\Delta_Y)+\frac{i}{2\pi}
\int_0^\lambda\Tr(S^\prime(s)S(-s))\;d\lambda$$
for $0\le\lambda<\sqrt{\mu_1}$.  Hence we get
$$\xi_c(0+)=-\frac{1}{4}(\Tr(S(0))+\dim\ker\Delta_Y).$$
Together with (\ref{2.3}) we obtain
$$\xi(0+)=-\dim\ker\Delta-\frac{1}{4}(\Tr(S(0))+\dim\ker\Delta_Y)$$
and by (\ref{2.2}) it follows that
$$b_0=\dim\ker\Delta+\frac{1}{4}(\Tr(S(0))+\dim\ker\Delta_Y).$$
Now observe that $S(0)$ satisfies $S(0)^2=\Id$. Hence
$$\Tr(S(0))+\dim\ker\Delta_Y=2\dim\ker(S(0)-\Id).$$ 
Combined with Proposition \ref{p2.1} we obtain the following corollary.
\begin{corollary}\label{c2.2}
Let $k=\dim\ker\Delta$ and $l=\dim\ker(S(0)-\Id)$. Then
$$\log\det(\Delta+\lambda,\Delta_0+\lambda)=(k+l/2)\log\lambda+\log\det(\Delta,
\Delta_0)+o(1)$$
as $\lambda\to0+$.
\end{corollary}

\section{Expansion of the Dirichlet-to-Neumann operator}
\setcounter{equation}{0}
Let $X=M\cup_Y Z$ be a manifold with a cylindrical end $Z=\R^+\times Y$ and
let $\Delta\colon C^\infty(Z,E)\to C^\infty(Z,E)$ be a Laplace type operator 
on $X$ with properties as above. For $z\in\C-\R_-$ let $R(z)$ be the 
Dirichlet-to-Neumann operator with respect to $\Delta+z$ and the
 hypersurface $Y=\{0\}\times Y\subset X$. 
In this section we study the expansion  of $\det(R(z))$ as $z\to0$. To begin
with we recall the definition of the Dirichlet-to-Neumann operator. Let
$z\in\C-R_-$ and $\varphi\in C^\infty(Y,E|Y)$. There exists a unique section
$\psi\in C^\infty(X-Y,E)\cap L^2(X,E)$ such that
\begin{equation*}
\begin{aligned}
(\Delta+z)\psi&=0\quad\text{on }\;X-Y;\\
\psi&=\varphi\quad \text{on }\; Y.
\end{aligned}
\end{equation*}
The solution $\psi$ is obtained as follows. Let $\widetilde\varphi\in 
C^\infty_c(X,E)$ be any extension of $\varphi$. Let $\Delta_D$ be the operator
$\Delta$ with Dirichlet boundary conditions along $Y$. Then
\begin{equation}\label{3.1a}
\psi=\widetilde\varphi-(\Delta_D+z)^{-1}
\bigl((\Delta+z)(\widetilde\varphi)\bigr).
\end{equation}
Furthermore, $\psi$ is continuous on $X$ and smooth on $\ov M$ and $\ov Z$. Its
normal derivative has a jump along $Y$. Then $R(z)\varphi$ is defined by
\begin{equation}\label{3.1}
R(z)\varphi=\frac{\partial}{\partial u}\bigl(\psi|_M\bigr)\big|_{\partial M}
-\frac{\partial}{\partial u}\bigl(\psi|_Z\bigr)\big|_{\partial Z}.
\end{equation}
By Theorem 2.1 of \cite{Ca}, $R(z)$ is an invertible pseudo-differential 
operator of order 1. Its principal symbol is given by
$$\sigma(R(z))(x,\xi)=2\sqrt{g_x(\xi,\xi)}\Id_{E_x},\quad \xi\in T_x^*Y.$$
Furthermore, $z\in\C-R_-\mapsto R(z)$ is a holomorphic function with values in
the space of pseudo-differential operators. Let $G(x,y,z)$ denote the kernel
of $(\Delta+z)^{-1}$. Then $G(x,y,z)$ is smooth in the complement of the
diagonal and for $x\not=y$, $G(x,y,z)\in\Hom(E_y,E_x)$. As shown in the proof of Theorem 2.1 of \cite{Ca}, we have
\begin{equation}\label{3.2}
R(z)^{-1}\varphi(x)=\int_Y G(x,y,z)\varphi(y)\;dy,\quad x\in Y,\; \varphi\in
C^\infty(Y,E|Y).
\end{equation}
In other words
$$R(z)^{-1}=\rho_Y\circ(\Delta+z)^{-1}(\cdot\otimes\delta_Y),$$
where $\rho_Y$ is the restriction map to $Y$ and $\delta_Y$ is the Dirac
 $\delta$-function along $Y$. Especially, if $\lambda>0$ then $R(\lambda)$ 
is an elliptic pseudodifferential operator of order 1 which is selfadjoint 
and positive definite. Hence its regularized determinant $\det( R(\lambda))$
 is defined.

Under the assumption that $\Delta_{M,D}$ is invertible, we can also define
the Dirichlet-to-Neumann operator with respect to $\Delta$ and $Y$. For this
purpose we need the following lemma.

\begin{lem}\label{l3.1}
For every $\varphi\in C^\infty(Y,E|Y)$ there exists a unique $\psi\in 
C^\infty(X-Y,E)\cap C^0(X,E)$, which is bounded and satisfies
\begin{equation}\label{3.3}
\begin{aligned}
\Delta\psi&=0\quad\text{on }\;X-Y;\\
\psi|_Y&=\varphi.
\end{aligned}
\end{equation}
\end{lem}
\begin{proof}
Since $\Delta_{M,D}$ is invertible, the Dirichlet problem on $M$ has a unique
solution, i.e., for every $\varphi\in C^\infty(Y,E|Y)$ there exists a unique
$\psi_1\in C^\infty(M,E)\cap C^0(\ov M,E)$ such that
\begin{equation*}
\begin{aligned}
\Delta_M\psi_1&=0\quad\text{in } \;M;\\
\psi_1|_Y&=\varphi.
\end{aligned}
\end{equation*}
Next we show that the Dirichlet problem on $Z$ has also a unique solution. Let
$\{\phi_i\}_{i\in\N}$ be an orthonormal basis of $L^2(Y,E|Y)$ consisting of 
eigenfunctions of $\Delta_Y$ with eigenvalues $0\le\lambda_0\le\lambda_1\le
\cdots$. Let $\varphi\in C^\infty(Y,E|Y)$. Then $\varphi$ has an expansion of
the form
$$\varphi=\sum_{i=1}^\infty a_j\phi_j.$$
Set 
$$\psi_2(u,y)=\sum^\infty_{j=1}a_je^{-u\sqrt{\lambda_j}}\phi_j(y).$$
Then $\psi_2\in C^\infty(Z,E)$ is bounded and satisfies
\begin{equation}\label{3.4}
\Delta\psi_2=0\;\mbox{and}\; \psi_2(0,y)=\varphi(y),\quad y\in Y.
\end{equation}
This proves existence. Now suppose that  $\widetilde\psi_2$ is a second
bounded solution of (\ref{3.3}). Set $g=\psi_2-\widetilde\psi_2.$ Then
$g\in C^\infty(Z,E)$ is bounded and satisfies
\begin{equation*}
\begin{split}
\left( -\frac{\partial^2}{\partial u^2}+\Delta_Y\right)g &  =0;  \\
g(u,y)  & =0,\;\;y\in Y.
\end{split}
\end{equation*}
If we expand $g$ in the orthonormal basis $\{\phi_j\}_{j\in\N}$ it 
follows that 
$$g(u,y)=\sum_{j=1}^m(b_ju+a_j)\phi_j(y)+\sum^\infty_{j=m+1}
(b_je^{{\sqrt{\lambda_j}}u}+a_je^{-\sqrt{\lambda_j}u})\phi_j(y),$$
where $m=\dim\ker\Delta_Y$. 
Since $g$ is bounded, it follows that $b_j=0$ for all $j\in\N$.
 Using that $g(0,y)=0,$ we obtain $a_j=0$ for all $j\in\N$. 
This proves uniqueness.
\end{proof}
Now we can proceed as above. Given $\varphi\in C^\infty(Y,E|Y),$ let
$\psi\in C^\infty(X-Y,E)\cap C^0(X,E)$ be the unique solution of 
(\ref{3.3}).
Then the Dirichlet-to-Neumann operator is defined by
\begin{equation}\label{3.5}
R\varphi=\frac{\partial}{\partial u}\bigl(\psi|_M\bigr)\big|_{\partial M}
- \frac{\partial}{\partial u}\bigl(\psi|_Z\bigr)\big|_{\partial Z}.
\end{equation}
Next we establish some properties of $R$.
\begin{lem}\label{l3.2}
There exist a smoothing operator $K$ such that
$$R=2\sqrt{\Delta_Y}+K.$$
\end{lem}
\begin{proof} Since $X-Y=M\sqcup Z$, $R$ can be written as
$$R=R_{\inte}+R_{\ext},$$
where $R_{\inte}$ is the Neumann jump operator on $M$. It is defined as
follows. Given  $\varphi\in C^\infty(Y,E|Y)$, let $\psi_1\in C^\infty(M,E)\cap
C^0(\ov M,E)$ be the unique solution of
\begin{equation*}
\Delta\psi_1=0 \quad \text{on }\;M,\quad \psi_1|_Y=\varphi.
\end{equation*}
Then $R_{\inte}$ is defined as 
$$R_{\inte}\varphi:=\frac{\partial\psi_1}{\partial u}\Big|_Y.$$
Similarly let $\psi_2\in C^\infty(Z,E)\cap C^0(\ov Z,E)$ be the unique
bounded solution of
$$\Delta\psi_2=0\quad\text{on }\;Z,\quad \psi_2|_Y=\varphi.$$
Set
$$R_{\ext}(\varphi):=-\frac{\partial\psi_2}{\partial u}\Big|_Y.$$
As explained above, $\psi_2$ is given by
\begin{equation*}
\psi_2(u,y)=\sum_{j=1}^\infty\langle\varphi,\phi_j\rangle 
e^{-\sqrt{\lambda_j}u}\phi_j(y).
\end{equation*}
Hence we get
$$-\frac{\partial\psi_2}{\partial u}(0,y)=\sum_{j=1}^\infty\sqrt{\lambda_j}
\langle\varphi,\phi_j\rangle\phi_j(y)=(\sqrt{\Delta_Y}\varphi)(y).$$
Thus $R_{\ext}=\sqrt{\Delta_Y}$. By Theorem 2.1 of \cite{Le3} it follows
that $R_{\inte}=\sqrt{\Delta_Y}+K$, where $K$ is a smoothing operator. 
This proves  the lemma.
\end{proof}
In particular, it follows that  $R$ is an elliptic  pseudodifferential
operator of order 1.
\begin{lem}\label{l3.2a}
For every $\phi\in C^\infty(Y,E|Y)$, $R(\lambda)\phi$ is a continuous function
of $\lambda\in [0,\infty)$ and
$$\lim_{\lambda\to 0+}R(\lambda)\phi=R\phi.$$
\end{lem}
\begin{proof}
Let $\lambda\ge 0$. As above, $R(\lambda)$ can be written as
$$R(\lambda)=R_{\inte}(\lambda)+R_{\ext}(\lambda).$$
Given $\phi\in C^\infty(Y,E|Y)$, let $\psi_1(\lambda)\in
C^\infty(M-Y,E)\cap C^0(M,E)$ be the unique section which satisfies 
$(\Delta+\lambda)\psi_1(\lambda)=0$ and $\psi_1(\lambda)|_Y=\phi$. Let
$\tilde\phi\in 
C^\infty(M,E)$ be any extension of $\phi$ which is smooth up to the boundary. 
Then
$$\psi_1(\lambda)=\tilde\phi-(\Delta_{M,D}+\lambda)^{-1}
((\Delta_M+\lambda)(\tilde\phi)).$$ 
Since $\Delta_{M,D}$ is invertible, this formula also holds for $\lambda=0$.
From this
representation of $\psi_1(\lambda)$ it follows immediately that
$R_{\inte}(\lambda)\phi$ converges to $R_{\inte}\phi$ as $\lambda\to 0+$.
Next observe that  the unique bounded solution $\psi_2(\lambda)\in 
C^\infty(Z,E)\cap C^0(Z,E)$ of 
$$(\Delta+\lambda)\psi_2(\lambda)=0\quad\mathrm{on}\;\;Z,
\quad\psi_2(\lambda)|_Y=\phi$$ 
is given by
\begin{equation*}
\psi_2(\lambda,u,y)=\sum_{j=1}^\infty\langle\varphi,\phi_j\rangle 
e^{-(\lambda_j+\lambda)^{1/2}u}\phi_j(y).
\end{equation*}
Then $R_{\ext}(\lambda)\phi:=\partial\psi_2(\lambda,u.y)/\partial u|_{u=0}$
and it follows that $R_{\ext}(\lambda)\phi$ is continuous in
$\lambda\in[0,\infty)$ 
and $R_{\ext}(\lambda)\phi$ converges to $R_{\ext}\phi$ as $\lambda\to 0+$.
\end{proof}

\begin{corollary}\label{c3.4}
The operator $R$ is formally selfadjoint and nonnegative.
\end{corollary}
\begin{proof}
As explained above, for every $\lambda>0$, the operator $R(\lambda)$ is
 formally selfadjoint and positive, and therefore the claim follows
 immediately from Lemma \ref{l3.2a}.
\end{proof}
Together we have proved that $R$ is a first order elliptic pseudo-differential
operator which is formally selfadjoint and nonnegative. Hence the regularized
determinant $\det R$ is well-defined.  

Our next purpose is  to study the bahaviour of the bounded operator
$R(\lambda)^{-1}$ 
as $\lambda\to0$. 
First we recall some facts about the spectral resolution of $\Delta$. For more
details we refer to \cite{Mu4}. We have
$$L^2(X,E)=L^2_d(X,E)\oplus L^2_c(X,E),$$
where
$$L^2_d(X,E)=\bigoplus_j\E(\lambda_j)$$
is the discrete sum of the eigenspaces of $\Delta$ with eigenvalues
$0\le\lambda_1<\lambda_1<\cdots.$ Each eigenspace is finite dimensional.
The orthogonal complement $L^2_c(X,E)$ of $L^2_d(X,E)$ is the absolutely 
continuous subspace for $\Delta$. It can be 
described in terms of generalized eigensections $E(\phi_j,\lambda)$
attached to the eigensections $\phi_j$ of $\Delta_Y.$ Each $E(\phi_j,\lambda)$ is a smooth section of $E$ and satisfies
$$\Delta E(\phi_j,\lambda)=\lambda E(\phi_j,\lambda).$$
Of particular importance for our purpose are the generalized eigensections
$E(\phi,\lambda)$ attached to $\phi\in\ker\Delta_Y.$ Let $\mu_1>0$ be the
smallest positive eigenvalue of $\Delta_Y.$ If we put $\lambda=s^2$ and
regard $E(\phi,\lambda)$ as a function of $s,$ then $E(\phi,s)$ has an 
analytic continuation to the disc $|s|<\mu_1.$
Let 
$$S(s):\ker\Delta_Y\to\ker\Delta_Y,\quad |s|<\mu_1,$$
be the corresponding scattering matrix. It is also holomorphic for 
$|s|<\mu_1$ and on $\R^+\times Y$ we have
\begin{equation}\label{3.6}
E(\phi,s,(u,y))=e^{ius}\phi(y)+e^{-isu}(S(s)\phi)(y)+\psi(s),
\end{equation}
where $\psi(s)$ is in $L^2.$ Let $0<\mu<\mu_1$ and let $P_\mu$ be the spectral
projection of $\Delta$ onto $[0,\mu]$.
By (\ref{3.2}) we have
\begin{equation}\label{3.7}
R(\lambda)^{-1}=\rho_Y\circ P_\mu(\Delta+\lambda)^{-1}(\cdot\otimes\delta_Y)+
\rho_Y\circ (\Id-P_\mu)(\Delta+\lambda)^{-1}(\cdot\otimes\delta_Y).
\end{equation}

First we study the second operator on the right. Let 
$$i_Y: L^2(Y,E|Y)\to H^{-1}(X,E)$$ be the map which is defined by
$i_Y(\varphi)=\varphi\delta_Y.$ 
Then $i_Y$ is continuous. Furthermore the restriction map $\rho_Y$ defines a 
continuous map 
$$\rho_Y\colon H^1(X,E)\to L^2(Y,E|Y).$$
Since $(\Delta+\lambda)^{-1} \colon H^{-1}(X,E)\to H^1(X,E)$ is continuous, we 
get a continuous map
$$\rho_Y\circ (\Id-P_\mu)(\Delta+\lambda)^{-1}\circ i_Y\colon L^2(Y,E|Y)\to 
L^2(Y,E|Y).$$

\begin{lem}\label{l3.3}
There exists $C>0$ such that
$$\parallel \rho_Y\circ (\Id-P_\mu)(\Delta+\lambda)^{-1}\circ
i_Y\parallel_{L^2} 
\le C$$
for all $\lambda\ge0$.
\end{lem}
\begin{proof}
Let $\varphi\in H^{-1}(X,E)$. Then $\parallel\varphi\parallel_{H^{-1}}=
\parallel(\Delta+\Id)^{-1/2}\varphi\parallel_{L^2}$. Hence we get
\begin{equation*}
\begin{split}
\parallel(\Id-P_\mu)(\Delta+\lambda)^{-1}\varphi\parallel_{H^1}&=
\parallel(\Delta+\Id)(\Id-P_\mu)(\Delta+\lambda)^{-1}
(\Delta+\Id)^{-1/2}\varphi\parallel_{L^2}\\
&\le \parallel(\Delta+\Id)(\Id-P_\mu)(\Delta+\lambda)^{-1}\parallel_{L^2}\cdot
\parallel\varphi\parallel_{H^{-1}}.
\end{split}
\end{equation*}
Using the spectral theorem we get
$$\parallel(\Delta+\Id)(\Id-P_\mu)(\Delta+\lambda)^{-1}\parallel_{L^2}
\le 1+1/\mu$$
for $\lambda\ge0$. This implies
$$\parallel(\Id-P_\mu)(\Delta+\lambda)^{-1}\parallel_{L(H^{-1},H^1)}\le
1+1/\mu$$
for $\lambda\ge0$.
Since $i_Y$ and $\rho_Y$ are continuous, the lemma follows.
\end{proof}

It remains to consider the first operator on the right hand side of 
(\ref{3.7}). This is a smoothing operator whose kernel $R(y_1,y_2,\lambda)$ 
can be
described as
follows. Let $\{\varphi_j\}$ be an
 orthonormal basis of eigensections of $\Delta$ with eigenvalues
$0\le\lambda_1\le\lambda_2\le\cdots$ and let $\phi_1,\ldots,\phi_m$ be an
orthonormal basis of $\ker\Delta_Y.$ Then it follows from the explicite
description of the spectral resolution of $\Delta$ (see \cite{Gu},
\cite{Mu4}) that
\begin{equation}\label{3.8}
\begin{split}
R(y_1,y_2,\lambda)=&\sum_{\lambda_j\le\mu}(\lambda_j+\lambda)^{-1}
\varphi_j(y_1)\otimes \varphi_j(y_2)\\
&+\frac{1}{2\pi}\sum^m_{j=1}\int^\mu_0(s^2+\lambda)^{-1}
E(\phi_j,s,y_1)\otimes E(\phi_j,-s,y_2)ds.
\end{split}
\end{equation}
We shall now determine the behaviour of this kernel as $\lambda\to 0$.
The behaviour of the first sum is obvious and we  only need to investigate the
second sum. 
\begin{lem}\label{l3.3a}
Let $\phi_1,...,\phi_m$ be an orthonormal basis of $\ker\Delta_Y$.
Then 
$$\sum_{j=1}^m E(\phi_j,s,y_1)\otimes E(\phi_j,-s,y_2)$$
is an even function of $s$, $|s|<\mu_1$.
\end{lem}
\begin{proof}
We recall that the generalized eigensections and the scattering
matrix satisfy the following functional equations. 
Let $\phi\in\ker\Delta_Y$. Then 
\begin{equation}\label{3.8a}
\begin{aligned}
E(\phi,-s)&=E(S(-s)\phi,s),\\
S(s)S(-s)&=\Id,\quad S(s)^t=S(s),\quad |s|<\mu_1.
\end{aligned}
\end{equation}
Let $\phi_1,...,\phi_m$ be an orthonormal basis of
 $\ker\Delta_Y$. Then there exist analytic functions $a_{ij}(s)$, 
$i,j=1,...,m$, defined in $|s|<\mu_1$, such that
\begin{equation}\label{3.8b}
S(s)\phi_i=\sum_{j=1}^m a_{ij}(s)\phi_j,\quad i=1,...,m.
\end{equation}
Using (\ref{3.8a}) and (\ref{3.8b}) we get
\begin{equation*}
\begin{aligned}
\sum_{j=1}^m E(\phi_j,-s,y_1)\otimes E(\phi_j&,s,y_2)=
\sum_{j=1}^m E(S(-s)\phi_j,s,y_1)\otimes E(S(s)\phi_j,-s,y_2)\\
&=\sum_{j=1}^m\sum_{k,l=1}^m a_{jk}(-s)a_{jl}(s)E(\phi_k,s,y_1)\otimes
E(\phi_l,-s,y_2).
\end{aligned}
\end{equation*}
By (\ref{3.8a}) the matrix $A(s)=(a_{ij}(s))_{i,j}$ is symmetric and
satisfies $A(-s)A(s)=\Id$. This implies
$$\sum_{j=1}^m E(\phi_j,-s,y_1)\otimes E(\phi_j,s,y_2)=
\sum_{j=1}^m E(\phi_j,s,y_1)\otimes E(\phi_j,-s,y_2)$$
as claimed.
\end{proof}
By Lemma \ref{l3.3a} there exists a smooth section $\widetilde E(s)$ of 
$E\boxtimes E$ over $X\times X$ which is holomorphic for $|s|<\mu$ such that
\begin{equation}\label{3.8c}
\begin{split}
\sum_{j=1}^m E(\phi_j,s,y_1)\otimes E(\phi_j,-s,y_2)&=
\sum_{j=1}^m E(\phi_j,0,y_1)\otimes E(\phi_j,0,y_2)\\
&\quad+s^2\widetilde E(s,(y_1,y_2)),\quad |s|<\mu.
\end{split}
\end{equation}
Note that
$$\int_0^{\mu}\frac{ds}{s^2+\lambda}=
\frac{\pi}{2\sqrt{\lambda}}-\frac{1}{\mu}+O(\lambda)$$
as $\lambda\to0$. Together with (\ref{3.8c}) we get
\begin{equation*}
\begin{split}
\frac{1}{2\pi}\sum_{j=1}^m\int^\mu_{0}(s^2+\lambda)^{-1}E(\phi_j,s,y_1)
&\otimes E(\phi_j,-s,y_2)ds\\
&=\frac{1}{4\sqrt{\lambda}} \sum_{j=1}^m
E(\phi,0,y_1)\otimes E(\phi_j,0,y_2)+O(1)
\end{split}
\end{equation*}
as $\lambda\to0$. 
To continue we consider the  the scattering matrix $S(0)$ at zero energy. 
It satisfies
$$S(0)^2=\Id.$$
Let $\phi\in\ker\Delta_Y$. If $S(0)\phi=\phi$ then it follows from (\ref{3.6})
that on $\R^+\times Y$ we have
$$E(\phi,0)=2\phi+\psi,$$
where $\psi\in L^2(\R^+\times Y,E).$
If  $S(0)\phi=-\phi,$ then $E(\phi,0)=0$ \cite[p. 209]{Mu2}. Let 
$$\ker\Delta_Y=V^+\oplus V^-$$
be the decomposition of $\ker\Delta_Y$ in the $\pm1$- eigenspaces of $S(0)$. 
Then $V^+$ equals the space of limiting values of extended solutions of 
$\Delta$ \cite{Mu2}. Let $\phi_1,\ldots,\phi_l$ be an orthonormal basis of 
$V^+$ and let $\varphi_1,...,\varphi_m$ be an orthonormal basis of 
$\ker\Delta$. Define the kernel $R_1$ by
\begin{equation}\label{3.9}
R_1(y_1,y_2,\lambda)=\frac{1}{\lambda}\sum_{j=1}^m\varphi_j(y_1)
\otimes\varphi_j(y_2)+
\frac{1}{4\sqrt{\lambda}}\sum_{j=1}^lE(\phi_j,0,y_1)\otimes 
E(\phi_j,0,y_2).
\end{equation}
Let $R_1(\lambda)\colon L^2(Y,E|Y)\to L^2(Y,E|Y)$ be the operator defined
by this kernel. Together with Lemma \ref{3.3} we obtain
\begin{prop}\label{p3.4}
There exists a bounded operator $R_2(\lambda):L^2(Y,E|Y)\to L^2(Y,E|Y)$ such 
that
$$R(\lambda)^{-1}=R_1(\lambda)+R_2(\lambda),\quad \lambda>0,$$
and $\parallel R_2(\lambda)\parallel$ is uniformly bounded as $\lambda\to0$.
\end{prop}
Let $\H\subset C^\infty(X,E)$ be the subspace spanned by $\ker\Delta$ and
$E(\phi_1,0),\ldots,E(\phi_l,0)$. Then $\H$ is the subspace of all bounded
 sections $\phi\in C^\infty(X,E)$ such that $\Delta\phi=0$. Set
$$\H_Y=\rho_Y(\H).$$
\begin{lem}\label{l3.5} The restriction map $\rho_Y:\H\to\H_Y$ is an isomorphism.
\end{lem}
\begin{proof} Let $\phi\in\H$. Then $\Delta\phi=0$ and $\phi$ is bounded. 
Suppose that $\rho_Y(\phi)=0$. This means that $\phi|_Y=0$. 
By the uniqueness of the 
Dirichlet problem, it follows that $\phi=0$. Thus $\rho_Y$ is injective and hence
an isomorphism.
\end{proof}

\begin{lem}\label{l3.6} $\ker R=\H_Y.$
\end{lem}
\begin{proof} Let $\varphi\in\H_Y.$ Then there exists $\psi\in\H$ with
$\psi|_Y=\varphi.$ Moreover $\psi$ is bounded and $\Delta
\psi=0.$ Thus $\psi$ is a solution of the Dirichlet problem (\ref{3.3}). Since
$\psi$ is smooth on $X,$ it follows that $R\varphi=0.$ Now suppose that
$\varphi\in\ker R.$ Then there exists a bounded solution $\psi$ of (\ref{3.3})
such that
$$\frac{\partial}{\partial u}\bigl(\psi|_M\bigr)\big|_{\partial M}
=\frac{\partial}{\partial u}\bigl(\psi|_Z\bigr)\big|_{\partial Z}.$$
This implies that $\Delta \psi=0$ in the sense of distributions. By elliptic
regularity it follows that $\psi\in C^\infty(X,E)$ and
$\Delta\psi=0.$ If we expand $\psi|_Z$ in the orthonormal basis
$\{\phi_j\}_{j\in\N}$ we get
$$\psi(u,y)=\sum^m_{j=1}a_j\phi_j(y)+\sum^\infty_{j=m+1}a_je^{-u\sqrt{\lambda_j}}\phi_j(y),$$
where $m=\dim\ker\Delta_Y$. Let $\phi=\sum^m_{j=1}a_j\phi_j.$ Then we get
$$\psi|_Z=\phi+\psi_1,$$
where $\psi_1\in L^2$. Put $\tilde\psi=\psi-\frac{1}{2}E(\phi,0)$. Then
it follows that $\tilde\psi\in\ker\Delta$.  This implies that $\psi\in\H.$
\end{proof}

Let $\langle\cdot,\cdot\rangle_Y$ be the inner product in $\H_Y$ induced by 
the inner product in $L^2(Y,E|Y)$. Let $\varphi_1,\ldots,\varphi_k$ be an
 orthonormal basis of $\ker\Delta$. Set $\psi_i=\rho_Y(\varphi_i)$, if $1\le
 i\le k$, and $\psi_{k+j}=\frac{1}{2}\rho_Y(E(\phi_j),0))$, if $1\le j\le l$.
Set $a_{ij}=\langle\psi_i,\psi_j\rangle_Y$, $1\le i,j\le k+l$ and let $A$
be the $(k+l)\times (k+l)$-matrix with entries $a_{ij}$, $i,j=1,...,k+l$.
Then the main result of this section is the following theorem.

\begin{theo}\label{th3.7} Let $k=\dim\ker\Delta$ and $l=\dim V^+.$ Then
$$\log\det R(\lambda)=(k+l/2)\log\lambda-\log\det A+\log\det R+O(\lambda)
$$
as $\lambda\to 0+$. 
\end{theo}
\begin{proof} The proof is analogous to the proof of Theorem B of \cite{Le1}.
Let
$$0\le\mu_1(\lambda)\le \cdots\le\mu_{k+l}(\lambda)<\mu_{k+l+1}
(\lambda)\le\cdots$$
be the eigenvalues of $R(\lambda).$ By Lemma \ref{l3.6} it follows that
$$\lim_{\lambda\to 0}\mu_j(\lambda)=0\quad\mbox{for}\;
1\le j\le k+l,$$
and $\mu_j(\lambda)\ge c>0$ for $j>k+l.$ Then
\begin{equation}\label{3.10}
\log\det R(\lambda)=\log(\mu_1(\lambda)\cdots\mu_{k+l}(\lambda))+
\log\det R+O(\lambda)
\end{equation}
as $\lambda\to0.$ So it remains to determine the behaviour of 
$\log(\mu_1(\lambda)\cdots\mu_{k+l}(\lambda))$ as $\lambda\to0.$
Let $\eta_1(\lambda),...,\eta_{k+l}(\lambda)$ be an orthonormal set of 
eigensections  of $R(\lambda)$ corresponding to the eigenvalues
$\mu_1(\lambda),...,\mu_{k+l}(\lambda)$. Let $1\le j\le k+l$.
By Proposition \ref{p3.4} we get
 \begin{equation*}
\mu_i(\lambda)^{-1}\delta_{ij}
=\langle R(\lambda)^{-1}\eta_i(\lambda),\eta_j(\lambda)\rangle
=\langle R_1(\lambda)\eta_i(\lambda),\eta_j(\lambda)\rangle
+\langle R_2(\lambda)\eta_i(\lambda),\eta_j(\lambda)\rangle,
\end{equation*}
and the second term on the right remains bounded as $\lambda\to0+$. 
By (\ref{3.9}) the first term equals
\begin{equation}\label{3.11}
\begin{aligned}
\langle R_1(\lambda)\eta_i(\lambda),\eta_j(\lambda)\rangle&=
\frac{1}{\lambda}\sum_{p=1}^k\langle\varphi_p,\eta_i(\lambda)\rangle_Y
\langle\varphi_p,\eta_j(\lambda)\rangle_Y\\
&\quad+\frac{1}{4\sqrt{\lambda}}\sum_{q=1}^l\langle E(\phi_q,0),\eta_i(\lambda)
\rangle_Y\langle E(\phi_q,0),\eta_j(\lambda)\rangle_Y.\\
\end{aligned}
\end{equation}
Set 
\begin{equation*}
\widetilde\psi_i(\lambda)=
\begin{cases}\rho_Y(\varphi_i),& \text{if }\; 1\le i\le k,\\
\frac{\lambda^{1/4}}{2}\rho_Y(E(\phi_{i-k},0)),& \text{if}\;
k+1\le i\le k+l. 
\end{cases}
\end{equation*}
Let $\widetilde a_{ij}(\lambda)=\langle\widetilde
\psi_i(\lambda),\eta_j(\lambda)\rangle$
and let $\widetilde A(\lambda)$ be the matrix with entries
 $\widetilde a_{ij}(\lambda)$, $1\le i,j\le k+l$. Then (\ref{3.11}) can be
written as
\begin{equation*}
\langle R_1(\lambda)\eta_i(\lambda),\eta_j(\lambda)\rangle=
\frac{1}{\lambda}(\widetilde A(\lambda)^t\widetilde A(\lambda))_{ij}
\end{equation*}
and we get
\begin{equation*}
\mu_i(\lambda)^{-1}\delta_{ij}=\frac{1}{\lambda}
(\widetilde A(\lambda)^t\widetilde A(\lambda))_{ij}+O(1)
\end{equation*}
 as $\lambda\to0+$. Note that 
$(\widetilde A(\lambda)^t\widetilde A(\lambda))_{ij}$ is bounded as 
$\lambda\to0+$. Hence for $i\not=j$ we get
$(\widetilde A(\lambda)^t\widetilde A(\lambda))_{ij}=O(\lambda)$ as 
$\lambda\to0+$. This implies
\begin{equation}\label{3.12}
(\mu_1(\lambda)\cdots\mu_{k+l}(\lambda))^{-1}=\lambda^{-(k+l)}\det
(\widetilde A(\lambda)^t\widetilde A(\lambda))(1+O(\lambda))
\end{equation}
as $\lambda\to0+$. Now observe that 
$\widetilde A(\lambda)\widetilde A(\lambda)^t$ is equal to the 
matrix with entries $\langle\widetilde\psi_i(\lambda),
\widetilde\psi_j(\lambda)\rangle$, $1\le i,j\le k+l$. Let
\begin{equation*}
C(\lambda)=\begin{pmatrix}\Id_k&0\\
0&\lambda^{1/4}\Id_l
\end{pmatrix}.
\end{equation*}
Then it follows from the definition of $A$ that
$$\widetilde A(\lambda)\widetilde A(\lambda)^t=C(\lambda)\cdot A\cdot 
C(\lambda).$$
Together with (\ref{3.12}) we obtain
\begin{equation*}
(\mu_1(\lambda)\cdots\mu_{k+l}(\lambda))^{-1}=\lambda^{-(k+l/2)}\det(A)
(1+O(\lambda)).
\end{equation*}
Taking the logarithm and inserting the result in (\ref{3.10}), the 
 theorem follows.

\end{proof}

\section{Proof of Theorem \ref{th1.1a}}
\setcounter{equation}{0}
Let $\lambda>0$. By Theorem 4.2 of 
\cite{Ca} there is a polynomial $P(\lambda)$ with real coefficients of 
degree $\le(n-1)/2$ such that
\begin{equation}\label{4.1}
\frac{\det(\Delta+\lambda,\Delta_0+\lambda)}{\det(\Delta_{M,D}+\lambda)}
=e^{P(\lambda)}\det (R(\lambda)).
\end{equation}
All terms have asymptotic expansions as $\lambda\to 0$. Since $\Delta_{M,D}$
is invertible, $\det(\Delta_{M,D}+\lambda)$ is continuous at $\lambda=0$ and
$\lim_{\lambda\to 0}\det(\Delta_{M,D}+\lambda)=\det(\Delta_{M,D})$.
Next consider the polynomial $P(\lambda)$. 
In the proof of
Proposition 4.7 of \cite{Ca}, Carron has shown that the polynomial 
$P(\lambda)$ can be computed in terms
of the coefficients of the asymptotic expansion of $\Tr(e^{-t\Delta_Y})$ as
$t\to0$. Let
$$\Tr(e^{-t\Delta_Y})\sim \sum_{j=0}^\infty a_j t^{-(n-1)/2+j}, \quad t\to0+,$$
be the heat expansion. If $n$ is even, we have $P=0$, and if $n=2p+1$ then
$$P(\lambda)=-\log(2) \sum_{j=0}^{p}\frac{(-1)^{p-j}}{(p-j)!}a_j
\lambda^{p-j}.$$
In particular, it follows that
$$P(0)=-\log(2)(h_Y+\zeta_Y(0)),$$
where $h_Y=\dim\ker\Delta_Y$ and $\zeta_Y(s)$ is the zeta function of 
$\Delta_Y$. Together with Corollary \ref{c2.2}  and Theorem \ref{th3.7}, 
Theorem \ref{th1.1a} follows.

\section{Regularized determinants on a finite cylinder}
\setcounter{equation}{0}

In this section we study the regularized determinant of a Laplace type 
operator on a finite cylinder over a closed Riemannian manifold $Y$. 
Let $\Delta_Y:C^\infty(Y,E_0)\to
C^\infty(Y,E_0)$ be a Laplace type operator on $Y$. For $r>0$ set
$$Z_r=[0,r]\times Y.$$
Let $E\to Z_r$ be the pull back bundle of $E_0$, i.e., 
$E=[0,r]\times E_0$. Let
$$\Delta=\Delta_{Z_r}:=-\frac{\partial^2}{\partial u^2}+\Delta_Y:
C^\infty(Z_r,E)\to C^\infty(Z_r,E).$$
Then $\Delta$ is a Laplace type operator on $Z_r$.
Impose Dirichlet boundary conditions at $\partial Z_r$ and let $\Delta_D$ be 
the corresponding self-adjoint extension. Then $\Delta_D$ is positive definite.
Let 
$$0\le\mu_1 \le\mu_2\le\cdots\to+\infty$$
be the eigenvalues of $\Delta_Y,$ counted with multiplicity. Let
$\zeta_Y(s)$ be the zeta function of $\Delta_Y$ and set
\begin{equation}\label{5.1}
\xi_Y(s)=\frac{\Gamma(s-1/2)}{\sqrt{\pi}\Gamma(s)}\zeta_Y(s-1/2).
\end{equation}
Sine $\zeta_Y(s)$ has at most a simple pole at $s=-1/2,$ $\xi_Y(s)$ is 
holomorphic at $s=0$. The main result of this section is the following 
proposition.

\begin{prop}\label{p5.1}
Let $h_Y=\dim\ker\Delta_Y.$
Then
\begin{equation}\label{5.2}
\det(\Delta_D)  = (2r)^{h_Y}e^{-r\xi^\prime_Y(0)/2}(\det \Delta_Y)^{-1/2}
\cdot\prod_{\mu_j>0}(1-e^{-2r\sqrt{\mu_j}}).
\end{equation}
\end{prop}
\begin{proof} The eigenvalues of $\Delta_D$ are given by
$$\lambda_{k,l}=\mu_l+\left(\frac{\pi}{r}\right)^2k^2,
\quad k,l\in\N.$$
Hence the zeta function of $\Delta_D$ equals
$$\zeta_{\Delta_D}(s)=\sum_{k,l\in\N}
\Bigl(\mu_l+\left(\frac{\pi}{r}\right)^2k^2\Bigr)
^{-s},\quad \Re(s)>\frac{d+1}{2},$$
where $d=\dim Y$. 
Let $\zeta(s)$ denote the Riemann zeta function. Then
\begin{equation}\label{5.3}
\zeta_{\Delta_D}(s) =h_Y\left(\frac{\pi}{r}\right)^{-2s}\zeta(2s)+
\sum_{k\in\N}\sum_{\mu_l>0}
\Bigl(\mu_l+\left(\frac{\pi}{r}\right)^2k^2\Bigr)^{-s},
\quad\Re(s)>\frac{d+1}{2}.
\end{equation}

Recall that $\zeta(0)=-1/2$ and $\zeta^\prime(0)=-1/2\log(2\pi)$. Hence
we get 
\begin{equation}\label{5.4}
\frac{d}{ds}\left\{\left(\frac{\pi}{r}\right)^{-2s}
\zeta(2s)\right\}\bigg|_{s=0}=
-\log2-\log r.
\end{equation}
Set
\begin{equation*}
\zeta_1(s):=\sum_{k\in\N}\sum_{\mu_l>0}
\Bigl(\mu_l+ \left(\frac{\pi}{r}\right)^2k^2\Bigr)^{-s},\quad 
\Re(s)>\frac{d+1}{2}.
\end{equation*}
By the Poisson summation formula we get
\begin{equation}\label{5.4a}
\begin{split}
\Gamma(s)\zeta_1(s) &= \sum_{\mu_l>0}\int^\infty_0e^{-\mu_lt}
\sum_{k\in\N} e^{-(\pi/r)^2k^2t}t^{s-1}dt\\
& = \sum_{\mu_l>0}\int^\infty_0e^{-\mu_lt}
\left( \frac{r}{\sqrt{\pi t}}\sum_{k\in\N}e^{-r^2k^2/t}+
\frac{1}{2}\left(\frac{r}{\sqrt{\pi t}}-1\right)\right)t^{s-1}dt\\
& = \frac{r}{\sqrt{\pi}}\sum_{k\in\N}\sum_{\mu_l>0}
\int^\infty_0 e^{-(\mu_lt+r^2k^2/t)}t^{s-3/2}dt\\
&\hskip2truecm+\frac{1}{2}\frac{r}{\sqrt{\pi}}\Gamma(s-1/2)\zeta_Y(s-1/2)
-\frac{1}{2}\Gamma(s)\zeta_Y(s).
\end{split}
\end{equation}

Denote by $T(s)$ the integral-series on the right hand side.
For $a,b,c \neq0$ and $s\in\C$ set
$K_s(a,b)=\int^\infty_0e^{-(a^2t+b^2/t)}t^{s-1}dt$ and
$K_s(c)=\int^\infty_0e^{-c(t+1/t)}t^{s-1}dt.$

It is proved in \cite[p.270f]{La} that the following relations hold
\begin{equation}\label{5.5}
K_s(c)=K_{-s}(c),\;  K_s(a,b)=\left(\frac{b}{a}\right)^sK_s(ab),
\;  K_{1/2}(c)=\sqrt\frac{\pi}{c}e^{-2c}.
\end{equation}

Furthermore, for every $x_0>0$ and $\sigma_0<\sigma_1$ there exists $C=C
(x_0,\sigma_0,\sigma_1)$ such that
\begin{equation}\label{5.6}
|K_s(x)|\le C e^{-2x}
\end{equation}
for all $x\ge x_0$ and $\Re(s)\in[\sigma_0,\sigma_1]$ \cite{La}.
With this notation we have
\begin{equation*}
T(s)=\frac{r}{\sqrt{\pi}}\sum_{k\in\N}\sum_{\mu_l>0}
K_{s-1/2}(\sqrt{\mu_l},rk).
\end{equation*}

Using (\ref{5.5}) and (\ref{5.6}) it follows that $T(s)$ is an entire
function of $s$.. Especially it is holomorphic at $s=0$. Since by (\ref{5.5}) 
we have $K_{-1/2}(a,b)=\frac{\sqrt{\pi}}{b}e^{-2ab},$ we get
\begin{equation}\label{5.7}
T(0)=\frac{r}{\sqrt{\pi}}\sum_{k\in\N}\sum_{\mu_l>0}
\frac{\sqrt{\pi}}{rk}e^{-2r\sqrt{\mu_l}k}=
-\sum_{\mu_l>0}\log(1-e^{-2r\sqrt{\mu_l}}).
\end{equation}
Thus by (\ref{5.4a}) we have
$$\zeta_1(s)=\frac{1}{\Gamma(s)}T(s)+\frac{1}{2}\frac{r}{\sqrt{\pi}}
\frac{\Gamma(s-1/2)}{\Gamma(s)}\zeta_Y(s-1/2)
-\frac{1}{2}\zeta_Y(s).$$
Using that $\xi_Y(s)$ is holomorphic at $s=0$, we obtain
\begin{equation*}
\zeta^\prime_1(0)=T(0)+r\xi^\prime_Y(0)-\frac{1}{2}\zeta^\prime_Y(0).
\end{equation*}
Together with (\ref{5.3}), (\ref{5.4}) and (\ref{5.7}) we get
\begin{equation*}
\zeta^\prime_{\Delta_D}(0)=-\sum_{\mu_k>0}\log\left(1-e^{-2r{\sqrt{\mu_k}}}
\right)
+r\xi^\prime_Y(0)-\frac{1}{2}\zeta^\prime_Y(0)-h_Y(\log 2+\log r).
\end{equation*}
This implies the claimed equality.
\end{proof}

\section{The decomposition formula}
\setcounter{equation}{0}

Let $(M,g)$ be a closed connected $n$-dimensional Riemannian manifold 
 and let $Y\subset M$ be a separating hypersurface
as in the introduction such that
$$M=M_1\cup_Y M_2,\quad M_1\cap M_2=Y.$$
We assume that the metric $g$ is a product on a tubular
neighborhood $N$ of $Y$. For $r\ge0$ let
$$M_{1,r}=M_1\cup([-r,0]\times Y),\quad
M_{2,r}=M_2\cup([0,r]\times Y),$$
where we identify $Y$ with $\{-r\}\times Y$ in the first case and with
$\{r\}\times Y$ in the second case. Set
$$M_r=M_{1,r}\cup_{\{0\}\times Y}M_{2,r},\quad N_r=[-r,r]\times Y.$$
Then
\begin{equation}\label{6.1a}
M_r=M_1\cup_Y N_r\cup_Y M_2,
\end{equation}
where $\partial M_1$ is identified with $\{-r\}\times Y$ and $\partial M_2$
with $\{r\}\times Y$.
The metric $g$ on $M$ has an obvious extension to a metric on $M_r$. 
Furthermore, let 
$$M_{i,\infty}=M_i\cup_Y(\R^+\times Y),\quad i=1,2.$$
Let $\Delta_M:C^\infty(M,E)\to C^\infty(M,E)$
be a Laplace type operator as in the introduction. and let $\Delta_{M_r}$ be
its 
canonical extension to a Laplace type operator on $M_r$, i.e. $\Delta_{M_r}$ 
is uniquely defined by
$$\Delta_{M_r}\big|_{M_i}=\Delta|_{M_i},\quad
\Delta_{M_r}\big|_{N_r}=-\frac{\partial^2}{\partial u^2}+\Delta_Y.$$
Let $\Delta_{M_i}=\Delta|_{M_i}$ and let $\Delta_{M_i,D}$ be
the selfadjoint extension of $\Delta_{M_i}:C^\infty_c(M_i,E)\to L^2(M_i,E)$
with respect to Dirichlet boundary conditions. We assume that $\Delta_{M_1,D}$ 
and $\Delta_{M_2,D}$ are invertible. Let $\Delta_{N_r,D}$
denote the selfadjoint extension of
$$-\frac{\partial^2}{\partial u^2}+\Delta_Y:C^\infty_c(N_r,E)\to L^2(N_r,E)$$
with respect to Dirichlet boundary conditions.
Let $Y_{\pm r}:=\{\pm r\}\times Y$ and denote by 
$\Sigma_r\subset M_r$ the hypersurface
$$\Sigma_r:=Y_{-r}\sqcup Y_r.$$
 Given $z\in\C-\R^-,$
let $R_{r}(z)$ be the Dirichlet-to-Neumann operator associated to 
 $(\Delta_{M_r}+z)$ 
and the hypersurface $\Sigma_r$. We recall the definition of 
$R_r(z)$. Let $\phi\in C^\infty(\Sigma_r,E_r|{\Sigma_r})$. There exists a
 unique section $\varphi\in C^\infty(M_r-\Sigma_r,E_r)\cap C^0(M_r,E_r)$ such
 that 
\begin{equation}\label{6.1c}
\begin{aligned}
(\Delta_{M_r}+z)\varphi&=0\quad \mathrm{on}\;\; M_r-\Sigma_r;\\
\varphi&=\phi\quad\mathrm{on}\;\;\Sigma_r.
\end{aligned}
\end{equation}
Then $R_r(z)(\phi)$ is given by
 \begin{equation}\label{6.1b}
\begin{aligned}
R_r(z)(\phi)\big|_{Y_{-r}}&=\frac{\partial}{\partial u}
\bigl(\varphi|_{M_1} \bigr)\big|_{\partial M_1}-
\frac{\partial}{\partial u}\bigl(\varphi|_{N_r}\bigr)\big|_{Y_{-r}},\\
R_r(z)(\phi)\big|_{Y_r}&=\frac{\partial}{\partial u}
\bigl(\varphi|_{N_r} \bigr)\big|_{Y_r}-
\frac{\partial}{\partial u}\bigl(\varphi|_{M_2}\bigr)\big|_{\partial M_2}.
\end{aligned}
\end{equation}

 Now we apply the Mayer-Vietoris
formula of \cite{BFK}, specialized to our case. We note that Theorem 1.4 of
 \cite{Ca} also holds in our case, where $M_r-\Sigma_r$ consists of three
 components. Thus
there exists a polynomial $P(z)$ with real coefficients of degree 
$<(n-1)/2$ such that for every $z\in\C-\R^-:$
$$\frac{\det (\Delta_{M_r}+z)}{\det(\Delta_{N_r,D}+z)\det(\Delta_{M_1,D}+z)
\det(\Delta_{M_2,D}+z)}=e^{P(z)}\det R_r(z).$$
Since we assume that the metric of $M_r$ is a product on a tubular 
neighborhood of $\Sigma_r,$ the polynomial depends only on $Y$ and can be 
computed as follows. Let $\zeta_Y(s,z)$ be the zeta function of $\Delta_Y+z.$ 
Then it follows from \cite[Theorem 6.3]{PW1} and also from the proof of
Proposition 4.7 of \cite{Ca} that
$$P(z)=-2\zeta_Y(0,z).$$
Thus 
\begin{equation}\label{6.1}
\frac{\det(\Delta_{M_r}+z)}{\det(\Delta_{N_r,D}+z)
\det(\Delta_{M_1,D}+z)\det(\Delta_{M_2,D}+z)}=2^{-2\zeta_Y(0,-z)}\det R_r(z).
\end{equation}

Now take $z=\lambda>0$ and consider the limit as $\lambda\to 0$ of the left
and right hand side of (\ref{6.1}).  Since  $\Delta_{M_i,D}$,
$i=1,2$,  and $\Delta_{N_r,D}$ are invertible,  it follows that 
\begin{equation}\label{6.2}
\lim_{z\to 0}\det(\Delta_{M_i,D}+\lambda)=\det\Delta_{M_i,D},\quad
\lim_{z\to0}\det(\Delta_{N_r,D}+\lambda)=\det\Delta_{N_r,D}.
\end{equation}
Let $h_r=\dim\ker\Delta_{M_r}.$ Then 
$$\det(\Delta_{M_r}+\lambda)=\lambda^{h_r}\det(\Delta_{M_r}
|_{(\ker\Delta_{M_r})^\perp}+\lambda)$$
and therefore we get 
\begin{equation}\label{6.4}
\lim_{\lambda\to 0}\det(\Delta_{M_r}+\lambda)\lambda^{-h_r}=\det\Delta_{M_r}.
\end{equation}
Also note that 
\begin{equation}\label{6.5}
\lim_{\lambda\to 0}\zeta_Y(0,\lambda)=\zeta_Y(0)+h_Y,
\end{equation}
where $h_Y=\dim\ker\Delta_Y.$ It remains to consider the limit of 
$\det R_r(\lambda)$ as $\lambda\to 0$.  Let 
$$\rho_r\colon C^\infty(M_r,E_r)\to C^\infty(\Sigma_r,E|{\Sigma_r})$$
denote the restriction operator.
Let ${\H}_{r}:=\rho_r(\ker\Delta_{M_r}).$
\begin{lem}\label{l6.1}
$$\rho_{r}:\ker\Delta_{M_r}\to{\H}_{r}$$
is an isomorphism.
\end{lem}

\begin{proof}
Let $\phi\in\ker\Delta_{M_r}$ and suppose that $\phi|_{\Sigma_r}=0.$ Let
$\psi=\phi|_{N_r}.$ Then $\Delta_{N_r}\psi=0$ and 
$\psi|_{\partial N_r}=0.$ Since $\Delta_{N_r,D}$ is invertible, it
follows that $\psi=0.$ In the same way we get $\phi|_{M_i}=0$, $i=1,2$, and 
hence $\phi=0$. Thus $\rho_{r}$ is injective and therefore an isomorphism.
\end{proof}

Let $\Delta_{M_r,D}$ be the selfadjoint extension of 
$$\Delta_{M_r}:C^\infty_c(M_r-\Sigma_r,E_r)\to L^2(M_r,E_r)$$ with respect to
Dirichlet boundary conditions. By our assumption, $\Delta_{M_r,D}$ is
invertible 
and hence, the Dirichlet-to-Neumann operator $R_r$ associated to  
$\Delta_{M_r}$ with
respect to $\Sigma_r\subset M_r$ can be defined in the same way as $R_r(z)$.
\begin{lem}\label{l6.2}
We have
$$\ker R_r=\rho_{r}(\ker\Delta_{M_r}).$$
\end{lem}
\begin{proof}
Let $\varphi\in\ker\Delta_{M_r}$ and let $\phi=\rho_r(\varphi)$. Then $\varphi$
is a solution of the Dirichlet problem  with boundary value $\phi$. 
Since $\varphi$ is smooth on $M_r$, it follows that $R_r(\phi)=0$. Now suppose 
that $\phi\in\ker R_r$. Then there exists $\varphi\in
C^\infty(M_r-\Sigma_r)\cap C^0(M_r)$ such that $\Delta_{M_r}\varphi=0$ on
$M_r-\Sigma_r$, $\varphi|_{\Sigma_r}=\phi$ and 
$$\frac{\partial}{\partial u}\bigl(\varphi|_{M_1})\big|_{\partial M_1}=
\frac{\partial}{\partial u}\bigl(\varphi|_{N_r})\big|_{Y_{-r}},\quad 
\frac{\partial}{\partial u}\bigl(\varphi|_{N_r})\big|_{Y_{r}}=
\frac{\partial}{\partial u}\bigl(\varphi|_{M_2})\big|_{\partial M_2}.$$
This implies that $\Delta_{M_r}\varphi=0$ in the distributional sense. By
elliptic regularity we conclude that $\varphi\in\ker\Delta_{M_r}$ and
$\rho_r(\varphi)=\phi$.
\end{proof}
Let $\varphi_1,...,\varphi_p$ be an orthonormal basis of
$\ker\Delta_{M_r}$. Set 
$$b_{ij}=\langle \rho_{r}(\varphi_i),\; \rho_{r}(\varphi_j)\rangle_{\Sigma_r}
,\quad i,j=1,...,p$$
and let
$$B_r=(b_{ij})^p_{i,j=1}.$$
Then $B_r$ is a symmetric invertible matrix.
\begin{prop}\label{p6.3} Let $h_r=\dim\ker\Delta_{M_r}$. 
Then
$$\log\det R_r(\lambda)=h_r\log\lambda-\log\det B_r+\log\det R_r+O(\lambda)$$
as $\lambda\to\infty$.
\end{prop}
\begin{proof} We use Lemma \ref{l6.2} and proceed in the same way as in
the proof of Theorem (\ref{th3.7}). 
\end{proof}
Combining (\ref{6.1})-(\ref{6.5}) and Proposition (\ref{p6.3}) we obtain
\begin{equation}\label{6.6}
\begin{split}
\log\det(\Delta_{M_r})= & \log\det\Delta_{N_r,D}+\log\det\Delta_{M_1,D}
+\log\det\Delta_{M_2,D}\\
& -\log\det B_r+\log\det R_r-2(\zeta_Y(0)+h_Y)\log 2.
\end{split}
\end{equation}

Let $Z=\R^+\times Y$ and let $\Delta_0$ be the selfadjoint extension of the
symmetric operator
$$-\frac{\partial^2}{\partial u^2}+\Delta_Y:C^\infty_c(Z,E)\to L^2(Z,E)$$
with respect to Dirichlet boundary conditions at $\partial Z=\{0\}\times Y.$
Let $R_{i,\infty}$ be the Dirichlet-to-Neumann operator for $\Delta_{i,\infty}$
with respect to the hypersurface $Y=\{0\}\times Y\subset M_{i,\infty}.$
Let $A_i$ be the Grahm matrix defined by the restrictions of the extended
$L^2$-solutions of $\Delta_{i,\infty}$ to $Y$ as in Theorem \ref{th3.7}.
By Theorem (\ref{th1.1a}) we have
\begin{equation*}
\begin{split}
\log\det(\Delta_{i,\infty},\Delta_0) & = \log\det R_{i,\infty}
+\log\det\Delta_{M_i,D}-\log\det A_i\\
& -(\zeta_Y(0)+h_Y)\log 2.
\end{split}
\end{equation*}
Together with (\ref{6.6}) we get
\begin{prop}\label{p6.4} Let the notation be as above. Then
\begin{equation}\label{6.7}
\begin{split}
\log\det\Delta_{M_r} =&\log\det\Delta_{N_r,D}-\log\det B_r+\log\det R_r\\
& +\sum^2_{i=1}(\log\det(\Delta_{i,\infty},\Delta_0)-\log\det
R_{i,\infty}+\log\det A_i).
\end{split}
\end{equation}
\end{prop}
Our next purpose is to study the behaviour of the various terms in this
equality as 
$r\to\infty.$ This, of course, will require additional assumptions. We begin
with the consideration of $\det R_r.$

To this end we need to describe the operator $R_r$ more explicitely. Let
$Q_i$ denote the Neumann jump operator on $M_i.$ In the proof of Lemma
(\ref{l3.2}) we established the following equality
\begin{equation}\label{6.8}
R_{i,\infty}=Q_i+\sqrt{\Delta_Y},\quad i=1,2.
\end{equation}
Let $P_0:L^2(Y,E|Y)$ denote the orthogonal projection onto
$\ker\Delta_Y.$ Let
$$h_r(x)=\frac{\sqrt{x}}{\sinh(2\sqrt{x}r)}.$$
Define
$$K_r:L^2(Y,E|Y)\oplus L^2(Y,E|Y)\to
L^2(Y,E|Y)\oplus L^2(Y,E|Y)$$
by
\begin{equation}\label{6.9}
K_r:=\left(\frac{1}{2r}P_0+h_r(\Delta_Y)P_0{^\perp}\right)
\begin{pmatrix} 
e^{-2\sqrt{\Delta_Y}r} & -{\Id} \\
 -{\Id} & e^{-2\sqrt{\Delta_Y}r}
\end{pmatrix}.
\end{equation}
Set
\begin{equation}\label{6.10}
R_\infty=\begin{pmatrix} R_{1,\infty} & 0\\
0 & R_{2,\infty}\end{pmatrix}.
\end{equation}
Recall that $\Sigma_r\cong Y\sqcup Y.$ Using (\ref{5.7}) and the formula at
the  
bottom of p. 4104 of \cite{Le3}, it follows that
$$R_r:C^\infty(Y,E|Y)\oplus C^\infty(Y,E|Y)\to
C^\infty(Y,E|Y)\oplus C^\infty(Y,E|Y)$$
is given by
\begin{equation}\label{6.11}
R_r=R_\infty +K_r.
\end{equation}
Next observe that $K_r$ is a trace class operator and its trace norm
$\parallel K_r\parallel_1$ satisfies
\begin{equation}\label{6.12}
\parallel K_r\parallel_1\underset{r\to\infty}{\longrightarrow}0.
\end{equation}
By Corollary \ref{c3.4}, $R_{i,\infty}$, $i=1,2$, 
 are selfadjoint nonnegative operators in $L^2(Y,E|Y)$.
\begin{lem}\label{l6.5} Suppose that $R_{i,\infty}>0,$ $i=1,2.$ Then
$$\lim_{r\to\infty}\det R_r=\det R_{1,\infty}\cdot\det R_{2,\infty}.$$
\end{lem}
\begin{proof} This is proved in \cite[Lemma 4.1 ]{Le3}. For the convenience of
  the 
  reader we recall the proof.
It follows from (\ref{6.10}) and the assumptions that 
$R_\infty>0$. By (\ref{6.12}) it follows that there exists $r_0>0$ such
that the operator
$R_\infty+tK_r$ is invertible for $0\le t\le 1$ and $r\ge r_0.$ Thus
\begin{equation*}
\begin{split}
\log\det(R_\infty+K_r)-\log \det R_\infty & = \int^1_0\frac{d}{dt}\log\det(R_\infty+tK_r)dt \\
& = \int^1_0\Tr((R_\infty+t K_r)^{-1}K_r)dt\le\frac{1}{2\lambda_0}
\parallel K_r\parallel_1,
\end{split}
\end{equation*}
where $\lambda_0>0$ is the smallest eigenvalue of $R_\infty.$
The lemma follows from (\ref{6.12}).
\end{proof}
Let ${\H}_i,$ $i=1,2,$ be the space of extended $L^2$-solutions of
$\Delta_{i,\infty}.$ By Lemma (\ref{l3.5}) and Lemma (\ref{l3.6}) it follows
that 
$R_{i,\infty}$ is invertible if and only if ${\H}_i=\{0\},$ and the latter 
condition is a consequence of $\ker\Delta_Y=\{0\}$ and 
$\ker\Delta_{i,\infty}=\{0\}$. Furthermore, if $R_{\infty}$ is invertible,
it follows from (\ref{6.11}) and (\ref{6.12}) that $R_r$ is invertible for
$r\ge r_0$. By Lemma \ref{l6.1} and Lemma \ref{l6.2}, $R_r$ is invertible
if and only if $\ker\Delta_{M_r}=\{0\}$.

Using  these observation together with Proposition (\ref{p6.4}) and 
Lemma (\ref{l6.5}), we obtain
\begin{corollary}\label{c6.6} Suppose that $\ker\Delta_Y=\{0\}$ and
  $\ker\Delta_{i,\infty}=0,$ $i=1,2$.  Then $\Delta_{M_r}$ is invertible for
$r\ge r_0$  and 
$$\lim_{r\to\infty}\frac{\det\Delta_{M_r}}{\det\Delta_{N_r,D}}=
\det(\Delta_{1,\infty},\Delta_0)\cdot\det(\Delta_{2,\infty},\Delta_0).$$
\end{corollary}

The asymptotic behaviour of $\det\Delta_{N_r,D}$ as $r\to\infty$ is described
by Proposition \ref{p5.1}. Using this result, Theorem \ref{th1.1} follows.

Next we consider a compact Riemannian manifold $(X_0,g)$ with a nonempty
boundary $Y.$ We assume that the metric is a product on a collar
neighborhood $N=(-\epsilon,0]\times Y$ of $Y$ in $X_0.$ Let
$$\Delta_{X_0}\colon C^\infty(X_0,E)\to C^\infty(X_0,E)$$
be Laplace type operator as above such that on $N$ it equals $-\partial
^2/\partial u^2+\Delta_Y.$ For $r>0$ set
$$Z_r=[0,r]\times Y,\quad\mathrm{and}\quad X_r=X_0\cup_Y Z_r,$$
where $\{0\}\times Y\subset Z_r$ is identified with $\partial X_0=Y$. 
Let $X_\infty=X_0\cup_Y(\R^+\times Y)$ be the corresponding manifold with 
a cylindrical end. We extend $\Delta_{X_0}$ in the obvious to Laplace type 
operators
$\Delta_{X_r}$ and $\Delta_\infty$ on $X_r$ and $X_\infty,$ respectively.
Let $\Delta_{X_r,D}$ and $\Delta_{Z_r,D}$ denote the  Dirichlet
Laplacians associated to $\Delta_{X_r}$ and $\Delta_{Z_r}$, respectively.
 Furthermore, let $R_r$ be the Dirichlet-to-Neumann operator
associated to the decomposition $X_r=X_0\cup Z_r$. Let $\lambda>0$. 
Then by Theorem 4.2 of \cite{Ca} we have 
\begin{equation}\label{6.15}
\frac{\det(\Delta_{X_r,D}+\lambda)}{\det(\Delta_{X_0,D}+\lambda)
\det(\Delta_{Z_r,D}+\lambda)}=2^{-\zeta_Y(0,\lambda)}\det R_r(\lambda).
\end{equation}
As above, we let $\lambda\to 0$. Note that $\Delta_{Z_r,D}$ is invertible.
Assume that $\Delta_{X_r,D}$ is invertible. Then as $\lambda\to0$,
the determinants converge to the determinants of $\Delta_{X_r}$ and 
$\Delta_{Z_r}$, respectively. Furthermore as in Lemma \ref{l6.2} it follows 
that 
\begin{equation}\label{6.15a}
\ker R_r=\rho_Y(\ker\Delta_{M_r,D}).
\end{equation}
Hence $R_r$ is also invertible and
$$\log\det R_r(\lambda)=\log\det R_r +O(\lambda)$$
as $\lambda\to0$. Thus taking the limit $\lambda\to0$ of both sides of 
(\ref{6.15}), we get
\begin{equation}\label{6.16}
\frac{\det \Delta_{X_r,D}}{\det \Delta_{X_0,D}
\det \Delta_{Z_r,D}}=2^{-\zeta_Y(0)}\det R_r.
\end{equation}
Now recall that by Theorem \ref{1.1a} we have
$$\log\det(\Delta_\infty,\Delta_0)=\log\det R_\infty+\log\det \Delta_{X_r,D}-
\log\det A- \log(2)\; (\zeta_Y(0)+h_Y).$$
Combining this equality with (\ref{6.16}) we obtain
\begin{equation}\label{6.17}
\begin{split}
\log\det \Delta_{X_r,D}=&\log\det\Delta_{N_r,D}+\log\det R_r 
-\log\det R_\infty\\
&+\det(\Delta_\infty,\Delta_0)+\log\det A.
\end{split}
\end{equation}
To study the behaviour of $\det R_r$ as $r\to\infty$, we proceed as above. Let
$$f_r(x)=\frac{\sqrt{x}}{\sinh(r\sqrt{x})},\quad x\in\R.$$
Define the operator
$$L_r\colon L^2(Y,E|Y)\to L^2(Y,E|Y)$$
by
\begin{equation}\label{6.17a}
L_r:=\left(\frac{1}{r}
  P_0+f_r(\Delta_Y)P_0^\perp\right)e^{-r\sqrt{\Delta_Y}},
\end{equation}
where $P_0$ denotes the orthogonal projection onto $\ker\Delta_Y$. Then
$$R_r=R_\infty+L_r.$$
Suppose that $\ker\Delta_Y=\{0\}$ and $\ker\Delta_\infty=\{0\}$. Then it
follows from Lemma \ref{l3.6} that $\ker R_\infty=\{0\}$ and 
 by Lemma 4.1 of \cite{Le3} we get 
$$\lim_{r\to\infty}\det R_r=\det R_\infty.$$
Since $\parallel L_r\parallel\to0$  as $r\to\infty$, it follows that $R_r$ is
 invertible for $r\ge r_0$. By (\ref{6.15a}) this 
implies that $\Delta_{M_r,D}$ is invertible for $r\ge r_0$.
Under the same assumptions we have $\det A=1$. Together with (\ref{6.17}) 
we get
\begin{prop}\label{p6.7}
Suppose that $\ker\Delta_Y=\{0\}$ and $\ker\Delta_\infty=\{0\}$. Then
$$\lim_{r\to\infty}\frac{\det\Delta_{X_r,D}}{\det\Delta_{Z_r,D}}
=\det(\Delta_\infty,\Delta_0).$$
\end{prop}
Using Proposition \ref{p5.1}, it follows that as $r\to\infty$, 
\begin{equation}\label{6.18}
\det \Delta_{X_r,D}\sim
e^{-r\xi_Y^\prime(0)/2}\left(\det\Delta_Y\right)^{-1/2}
\det(\Delta_\infty,\Delta_0).
\end{equation}
We apply (\ref{6.18})  to $\det\Delta_{M_{i,r},D}$, $i=1,2$, and compare the
asymptotic behaviour with (\ref{1.6}). In this way we get
$$\lim_{r\to\infty}\frac{\det\Delta_{M_r}}{\det\Delta_{M_{1,r},D}
\det\Delta_{M_{2,r},D}}=\left(\det\Delta_Y\right)^{1/2}$$
which is the statement of Corollary \ref{c1.5}.

\section{Bochner-Laplace operators}
\setcounter{equation}{0}
In this section we study the case where $\Delta$ is a connection Laplacian.
To begin with we consider a manifold with a cylindrical end
$X=M\cup_Y Z$, $Z=\R^+\times Y$.
Let $F\to X$ be a Hermitian vector bundle over $X$ such that
$F|_{\R^+\times Y}=\mathrm{pr}_Y^*(F_0)$ for some Hermitian vector
bundle $F_0$ over $Y$. Let $\nabla$ be a metric connection in $F$ such that
on $\R^+\times Y$ it has the form
\begin{equation}\label{7.1}
\nabla=d_u\otimes \Id+d\otimes\nabla^Y,
\end{equation}
where $\nabla^Y$ is a metric connection on $F_0.$ Let
$$\Delta=\nabla^*\nabla,\quad  \Delta_Y=(\nabla^Y)^*\nabla^Y$$
be the associated Bochner-Laplace operators. Then
\begin{equation}\label{7.2}
\Delta=-\frac{\partial^2}{\partial u^2}+\Delta_Y\,\quad\mbox{on}\;\;
\R^+\times Y.
\end{equation}
Let $\overline\Delta$ be the unique selfadjoint extension of $\Delta_X$ in 
$L^2$.

Let $\{\phi_i\}_{i\in\N}$ be an orthonormal basis of $L^2(Y,F|Y)$ consisting
of eigensections of $\Delta_Y$ with eigenvalues
$$0\le\mu_1\le\mu_2\le\cdots\to+\infty.$$
\begin{lem}\label{l7.1}
We have
$$\ker\overline\Delta=\{0\}.$$
\end{lem}
\begin{proof} Let $\varphi\in C^\infty(X,F)$ be a square integrable solution
of $\Delta\varphi=0.$ Then
$$0=\langle\nabla^*\nabla\varphi,\varphi\rangle=\parallel\nabla \varphi\parallel^2.$$
Thus $\nabla\varphi=0.$ Since $\varphi$ is square integrable, it has the 
following expansion on $\R^+\times Y$ in terms of the orthonormal basis
$\{\varphi_i\}_{i\in\N}:$
\begin{equation}\label{7.3}
\varphi(u,y)=\sum_{\mu_i>0}c_ie^{-\sqrt{\mu_i}u}\phi_i(y),\quad
u\in\R^+,y\in Y.
\end{equation}
Furthermore,
$$\frac{\partial\varphi}{\partial u}(u,y)=\nabla_{\frac{\partial}{\partial u}}
\varphi=0.$$
Using (\ref{7.3}), it follows that the restriction of $\varphi$ to 
$\R^+\times Y$ vanishes. Since $\nabla\varphi=0$ and $\nabla$ is a metric
connection, it follows that $d\parallel\varphi\parallel^2=0$ and hence 
$\varphi= 0$.
\end{proof} 
Let
$$\ker\Delta_Y=V^+\oplus V^-$$
be the decomposition into the $\pm1$-eigenspaces of the scattering matrix
$S(0)$ (cf. $\S2$) and let
$$R:C^\infty(Y,F|Y)\to C^\infty(Y,F|Y)$$
be the Dirichlet-to-Neumann operator with respect to the hypersurface
$$Y=\{0\}\times Y\subset X.$$
\begin{lem}\label{l7.2}
We have
$$\ker R=V^+.$$
\end{lem}
\begin{proof}
Let $\H\subset C^\infty(X,E)$ be the space of bounded solutions of 
$\Delta\varphi=0$. By Lemma (\ref{l3.6}) we have
$\ker R=\rho_Y(\H).$
So it suffices to prove that 
$$\rho_Y(\H)=V^+.$$
Let $\varphi\in\H$. For $r>0$ let $X_r=M\cup_Y([0,r]\times Y)$. 
Using integration by parts, we get
\begin{equation}\label{7.4}
\begin{split}
0 & = \int_{X_r}\langle\nabla^*\nabla\varphi(x),\varphi(x)\rangle dx\\
& = \int_{X_r}|\nabla\varphi(x)|^2dx-
\int_Y\biggl\langle\frac{\partial}{\partial u}
\varphi(u,y),\varphi(u,y)\biggr\rangle\bigg|_{u=r} dy.
\end{split}
\end{equation}

Since $\varphi$ is  bounded and satisfies $\Delta\varphi=0$, it has  
the following  expansion on  $\R^+\times Y$:
\begin{equation}\label{7.5}
\varphi(u,y)=\sum^{h_Y}_{i=1}a_i\phi_i(y)
+\sum^\infty_{i=h_Y+1}b_ie^{-\sqrt{\mu_i}u}\phi_i(y), 
\end{equation}
where $h_Y=\dim\ker\Delta_Y$. This implies that the second integral on the
right of (\ref{7.4}) is exponentially decreasing as $r\to\infty$. Hence
$\nabla\varphi=0$. In particular, it follows that
$$\frac{\partial}{\partial u}\varphi(u,y)=0,\quad u\in \R^+,\; y\in Y.$$
Together with (\ref{7.5}) we get
\begin{equation}\label{7.6}
\varphi|_Z=\phi\in\ker\Delta_Y.
\end{equation}
Thus $\rho_Y(\H)\subset \ker\Delta_Y.$
Now recall that $\varphi\in\H$ if and only if there exist $\phi\in V^+$
and $\psi\in L^2(Z,F)$ such that $\varphi|_Z=\phi+\psi.$ By 
(\ref{7.6}) it follows that $\psi=0$. This proves that
$\rho_Y(\H)=V^+$.
\end{proof}
Let $A$ be the matrix that occurs in Theorem \ref{th1.1a}.
\begin{corollary}\label{c7.3} We have
$$ \det A=1.$$
\end{corollary}
\begin{proof} Recall the definition of $A.$ Given $\phi\in V^+,$ let  
$\frac{1}{2} E(\phi,0)$ be the extended solution of $\Delta_X$ with limiting
value 
$\phi$. Let $\phi_1,\ldots,\phi_p$ be an orthonormal basis of $V^+$. Let
$\psi_j=\frac{1}{2} \rho_Y(E(\phi,0))$.
Since by Lemma (\ref{l7.1}) $\ker\overline\Delta =\{0\},$ it follows that the
entries of $A$ are $a_{ij}=\langle\psi_i,\psi_j\rangle_Y.$ By
Lemma (\ref{l7.2}), we have $\frac{1}{2} \rho_Y(E(\phi,0))=\phi$ for $\phi\in
V^+$. Hence $a_{ij}=\delta_{ij}$.
\end{proof}

Now consider a compact Riemannian manifold $(M,g)$ and a Hermitian vector
bundle $E\to M$ as in the previous section. Let $\nabla$ be a metric
connection on $E$ such that on the tubular neighborhood $N=[-1,1]\times Y$ of
$Y$ in $M$
\begin{equation}\label{7.8}
\nabla=d_u\otimes \Id+\Id\otimes\nabla^Y,
\end{equation}
where $\nabla^Y$ is a metric connection on $E_0=E|Y.$ Let
$$\Delta_M=\nabla^*\nabla.$$
Let $E_r\to M_r$ and $E_{i,\infty}\to M_{i,\infty}$ be the canonical extensions
of the vector bundle $E\to M$ to vector bundles over $M_r$ and $M_{i,\infty},$ 
respectively. By (\ref{7.8}), $\nabla$ has a canonical extension to a 
connection $\nabla^r$ on $E_r$ and $\nabla^{i,\infty}$ on $E_{i,\infty,}$
$i=1,2,$ respectively.
Then 
$$\Delta_{M_r}=(\nabla^r)^*\nabla^r,\quad \Delta_{i,\infty}=(\nabla^{i,\infty})
^*\nabla^{i,\infty},\quad i=1,2.$$
Recall that the Dirichlet-to-Neumann operator $R_r$ is a selfadjoint operator
in $C^\infty(Y,E|Y)\oplus C^\infty(Y,E|Y)$.

Next we determine $\ker R_r$. Let $V_i^+\subset \ker\Delta_Y$ be the subspace
of limiting values of extended $L^2$-sections of $\Delta_{i,\infty}$, $i=1,2$. 
\begin{lem}\label{l7.4} We have
$$\ker R_r=\{(\phi,\phi)\mid \phi\in V_1^+\cap V_2^+\}.$$
\end{lem}
\begin{proof} By Lemma \ref{l6.2} we have
$$\ker R_r=\rho_{r}(\ker\Delta_{M_r}).$$
Let $\varphi\in\ker\Delta_{M_r}.$ Then
$$0=\langle\nabla^*\nabla\varphi,\varphi\rangle
=\parallel\nabla\varphi\parallel^2.$$
Thus
\begin{equation}\label{7.9}
\nabla\varphi=0 \quad \mbox{for all}\,\; \varphi\in\ker\Delta_{M_r}.
\end{equation}
Next observe that the restriction of $\varphi$ to $N_r$ satisfies
$$\left(-\frac{\partial^2}{\partial u^2}+\Delta_Y\right)\varphi(u,y)=0,
\quad u\in[-r,r],\; y\in Y,$$
and hence the expansion of $\varphi|_{N_r}$ in the orthonormal
basis $\{\phi_i\}_{i\in\N}$ is of the form
$$\varphi(u,y)=\sum^{h_Y}_{i=1}(a_i+b_i u)\phi_i(y)+\sum^{\infty}_{i=h_Y+1}
(c_i e^{-\sqrt{\mu_i}u}+d_i e^{\sqrt{\mu_i}u})\phi_i(y).$$
By (\ref{7.9}) it follows that $\frac{\partial}{\partial u}\varphi(u,y)=0,$ 
$u\in[-r,r].$ Hence we get
\begin{equation}\label{7.10}
\varphi(u,y)=\sum^{h_Y}_{i=1}a_i\phi_i(y),\quad (u,y)\in N_r.
\end{equation}
Actually, by our assumptions this holds on a slightly larger collar
neighborhood  of $Y$.
Denote the right hand side by $\phi$. Then $\phi\in\ker\Delta_Y$ and it
follows that
$$\rho_{r}(\varphi)=(\varphi(-r,\cdot),\;\varphi(r,\cdot))=(\phi,\phi).$$
Furthermore, let $\varphi_i=\varphi|M_i$, $i=1,2$. Then $\varphi_i$ satisfies
$$\Delta_{M_i}\varphi_i=0\quad\mathrm{and}\quad\varphi_i|_{(-\varepsilon,0]
\times Y}=\phi.$$
Thus $\varphi_i$, $i=1,2$, has a unique extension $\hat\varphi_i$ to an
extended $L^2$-solution of $\Delta_{i,\infty}$ with limiting value $\phi$.
This implies $\phi\in V_1^+\cap V_2^+$. On the other hand, suppose that
$\phi\in V_1^+\cap V_2^+$. Let $\hat\varphi_i\in C^\infty(M_{i,\infty},
E_{i,\infty})$, $i=1,2$,  be an extended $L^2$-solution with limiting value
 $\phi$.  By (\ref{7.6}) we have
$$\hat\varphi_i|_{(-\varepsilon,0]\times Y}=\phi.$$
Thus we can patch together $\hat\varphi_1|M_1$ and $\hat\varphi_2|M_2$
to a section $\varphi\in\ker\Delta_{M_r}$ with
$\rho_r(\varphi)=(\phi,\phi).$
\end{proof}

\begin{lem}\label{l7.4a} For all $r>0$ there exists an isomorphism
$$j_r:\ker\Delta_{M_r}\to\ker\Delta_M.$$
\end{lem}
\begin{proof}
Let $\varphi\in\ker\Delta_{M_r}.$ By (\ref{7.10}), there exists
$\phi\in\ker\Delta_Y$ such that 
\begin{equation}\label{7.11}
\varphi(u,y)=\phi(y),\quad (u,y)\in N_r.
\end{equation}
Note that by our assumption, $\nabla$ is the product connection on a slightly
larger tubular neighborhood $N_{r+\epsilon}$ of $Y$ and (\ref{7.11}) continues
to hold on $N_{r+\epsilon}.$ Set
$$\psi_i=\varphi|_{M_i},\quad i=1,2.$$
By (\ref{7.11}), it follows that 
\begin{equation}\label{7.12}
\psi_1\big|_{\partial M_1}=\psi_2\big|_{\partial M_2}.
\end{equation}
Define $\psi\in C^\infty(M-Y,E)\cap C^0(M,E)$
by
\begin{equation*}
\psi(x)=\begin{cases} \psi_1(x), & x\in M_1, \\
\psi_2(x), & x\in M_2.
\end{cases}
\end{equation*}
By the above observation, there exists a tubular neighborhood
$N_\epsilon=(-\epsilon,\epsilon)\times Y$ of $Y$ such that
$\psi|_{ N_\epsilon}=\phi$. Hence $\psi\in C^\infty(M,E)$ and
$\Delta_M\psi=0.$ By construction, the map
$$j_r:\varphi\in\ker\Delta_{M_r}\longmapsto\psi \in\ker\Delta_M$$
is injective and the inverse map can be defined in the same way. This proves
that $j_r$ is surjective.
\end{proof}
\begin{corollary}\label{c7.5} The dimension of $\ker_{M_r}$ is independent
of $r$. 
\end{corollary}

Put
$$q:=\dim\ker\Delta_{M_r}.$$
Let the matrix $B_r$  be defined as in the previous section.
Our next purpose is to study the behaviour of $\det B_r$ as $r\to\infty.$
To this end we need some auxiliary result. Let
$$P_i:C^\infty(Y,E|Y)\to C^\infty(M_i,E),\quad i=1,2$$
be the Poisson operator. Recall that for $\phi\in C^\infty(Y,E|Y)$,
$P_i(\phi)$ is the unique solution of the Dirichlet problem 
$$\Delta_{M_i}\psi=0,\quad \psi|_{\partial M_i}=\phi.$$

\begin{lem}\label{l7.6} There exists $C>0$ such that
$$\parallel P_i(\phi)\parallel\le C\parallel\phi\parallel,\quad
\phi\in\ker\Delta_Y,\quad i=1,2.$$
\end{lem}

\begin{proof}
There exists a collar neighborhood $(-\epsilon,0]\times Y$ of $Y$ in $M_i$ 
such that
\begin{equation}\label{7.13}
\Delta_{M_i}=-\frac{\partial^2}{\partial u^2}+\Delta_Y\quad\mbox{on}\;
(-\epsilon,0]\times Y.
\end{equation}
Let $f\in C^\infty(\R)$ be such that $f(u)=1$ for $u\ge -\epsilon/4$ and
$f(u)=0$ 
for $u\le-\epsilon/2$. Given $\phi\in C^\infty(Y,E|Y)$, set 
$$\widetilde\phi(u,y)=f(u) \phi(y),\quad u\in(-\epsilon,0],\;y\in Y,$$
and extend $\widetilde\phi$ by zero to a smooth section of $E\to M_i.$ Then
\begin{equation}\label{7.14}
P_i(\phi)=\widetilde\phi-(\Delta_{M_i,D})^{-1}(\Delta_{M_i}\widetilde\phi).
\end{equation}
Let $\phi\in\ker\Delta_Y.$ By (\ref{7.13}) we get
$$\Delta_{M_i}\widetilde\phi=-g^{\prime\prime}\phi.$$
Let $\lambda_1>0$ be the smallest eigenvalue of $\Delta_{M_i,D}$. Then by
(\ref{7.14}) we get
$$\parallel P_i(\phi)\parallel\le C_1\parallel\phi\parallel+\frac{1}{\lambda_1}
C_2\parallel\phi\parallel\le C\parallel\phi\parallel,$$
where $C>0$ is independent of $\phi\in\ker\Delta_Y.$
\end{proof}
\begin{lem}\label{l7.7}
Let $q=\dim\ker R_r.$ Then
$$r^q\det B_r=1+O(r^{-1})$$
as $r\to\infty$.
\end{lem}
\begin{proof} Let $\psi_{r,1},\ldots,\psi_{r,q}\in\ker\Delta_{M_r}$ be an
orthonormal basis of $\ker\Delta_{M_r}.$ Then $B_r$ is defined as
$$B_r=\left(\langle \rho_{r}(\psi_{r,i}),\rho_{r}(\psi_{r,k})\rangle\right)^
q_{i,k=1}.$$
By (\ref{7.10}), for each $r>0$ and
$k,k=1,\ldots,q,$ there exists $\phi_{r,k}\in\ker\Delta_Y$ such that
\begin{equation}\label{7.15} 
\psi_{r,k}(u,y)=\phi_{r,k}(y),\quad u\in[-r,r],\;y\in Y.
\end{equation}
Let $M_0=M_1\sqcup M_2.$ Then
\begin{equation}\label{7.16}
\delta_{ik}=\langle\psi_{r,i},\psi_{r,k}\rangle_{M_r}=\bigl\langle\psi_{r,i}
\big|_{M_0},\psi_{r,k}\big|_{M_0}\bigr\rangle_{M_0}+2r\langle
\phi_{r,i},\phi_{r,k}\rangle.
\end{equation}
By (\ref{7.15}) we have
$$\rho_{r}(\psi_{r,k})=(\phi_{k,r},\phi_{k,r}).$$
Hence by (\ref{7.16}) we get
\begin{equation}\label{7.17}
\langle \rho_{r}(\psi_{r,i}),\rho_{r}(\psi_{r,k})\rangle=2\langle\phi_{r,i},
\phi_{r,k}\rangle=\frac{1}{r}\left(1-\bigl\langle\psi_{r,i}\big|_{M_0},
\psi_{r,k}\big|_{M_0}\bigr\rangle\right).
\end{equation}
Furthermore, by (\ref{7.16})
\begin{equation}\label{7.18}
\parallel\phi_{r,k}\parallel^2\le\frac{1}{2r}.
\end{equation}
Now observe that by (\ref{7.15}) we have
$$\psi_{k,r}\big|_{\partial M_i}=\phi_{k,r},\quad k=1,\ldots,q.$$
Moreover $\Delta_{M_i}\psi_{k,r}=0$.
 Thus
$$\psi_{r,k}\big|_{M_0}=\psi_{r,k}\big|_{M_1}
+\psi_{r,k}\big|_{M_2}=P_1(\phi_{k,r})+P_2(\phi_{k,r}).$$ 
Together with Lemma (\ref{l7.6}) and (\ref{7.18}) it follows that there exists
$C>0$ such that
$$\parallel\psi_{k,r}\big|_{M_0}\parallel\le
C\parallel\phi_{k,r}\parallel\le\frac{c}{2\sqrt{r}}$$ 
for all $r>0$ and $k=1,\ldots,q$. 
Hence by (\ref{7.16}) we get
$$\langle \rho_{r}(\psi_{r,i}),\rho_{r}(\psi_{r,k})\rangle=\frac{1}{r}
\left(\delta_{ik}+O\left(\frac{1}{r}\right)\right)$$
as $r\to\infty.$ This implies
$r^q\det B_r=1+O(r^{-1}).$
\end{proof}

Our next purpose is to study the behaviour of $\det R_r$ as
$r\to\infty$. Recall that by Lemma \ref{l7.2}
\begin{equation}\label{7.19}
\ker R_\infty=V_1^+\oplus V_2^+.
\end{equation}
Furthermore, by Lemma \ref{l7.4} we have 
\begin{equation}\label{7.20}
\ker R_r=\{(\phi,\phi)\mid \phi\in V_1^+\cap V_2^+\}.
\end{equation}
To study $R_r$ on the orthogonal complement of $\ker R_r$ we need to introduce
some auxiliary subspaces of $L^2(Y,E|Y)\oplus L^2(Y,E|Y)$. First put
\begin{equation}\label{7.21}
L=(V_1^+\cap V_2^+)\oplus(V_1^+\cap V_2^+).
\end{equation}
By (\ref{7.19}) we have $L\subset \ker R_\infty$. Furthermore, it follows from
(\ref{6.9}) that on $\ker\Delta_Y\oplus\ker\Delta_Y$ the operator $K_r$ is
given by 
\begin{equation}\label{7.22}
K_r=\frac{1}{2r}\begin{pmatrix} \Id & -\Id\\
-\Id & \Id\end{pmatrix}.
\end{equation}
This implies that $L$ is invariant under $K_r$. Therefore, $L$ is an invariant
subspace for $R_r=R_\infty+K_r$. Let
$$W=\{(\phi,-\phi)\mid\phi\in V_1^+\cap V_2^+\}.$$
Then by (\ref{7.20}) we get an orthogonal decomposition
$$L=\ker R_r\oplus W$$
and it follows from (\ref{7.21}) that $W$ is an invariant subspace of $K_r$
and hence of $R_r$. Moreover
$$R_r|_W=\frac{1}{r}\Id.$$
Set
$$h_{12}:=\dim(V_1^+\cap V_2^+).$$
Note that $h_{12}=q=\dim \ker R_r$. Let $L^\perp$ be the orthogonal complement
 of $L$ in $L^2(Y,E|Y)\oplus 
L^2(Y,E|Y)$.
Then it follows that
\begin{equation}\label{7.22a}
\det R_r=r^{-h_{12}}\det\left(R_r|L^\perp\right).
\end{equation}
So we can continue with the investigation of $R_r|L^\perp$. Let $L_1\subset 
V_1^+\oplus V_2^+$ be the orthogonal complement of $L$ in $V_1^+\oplus V_2^+$
and $(\ker R_\infty)^\perp$ the orthogonal complement of $\ker R_\infty$ in
$L^2(Y,E|Y)\oplus L^2(Y,E|Y)$.
Then
\begin{equation}\label{7.23}
L^\perp=L_1\oplus(\ker R_\infty)^\perp
\end{equation}
with $L_1\subset \ker R_\infty$. This decomposition is invariant under
 $R_\infty$, however, it is not invariant under $K_r$ and hence, it is not
 invariant under $R_r$. In fact, with respect to (\ref{7.23}) we may write
$$R_r|_{L^\perp}=\begin{pmatrix} A(r) & B(r)\\
C(r) & D(r)\end{pmatrix},$$
where the operators $A(r),...,D(r)$ are defined as follows. Let $\Pi_1$ denote
the orthogonal projection of $L^\perp$ onto $L_1$. Then
$$A(r)=\Pi_1 K_r\Pi_1,\quad B(r)=\Pi_1 K_r(\Id-\Pi_1),\quad
C(r)=(\Id-\Pi_1) K_r\Pi_1,$$
and 
\begin{equation}\label{7.23a}
D(r)=R_\infty|_{(\ker R_\infty)^\perp}+(\Id-\Pi_1) K_r(\Id-\Pi_1).
\end{equation}
Recall that $K_r$ is a trace class operator whose trace norm $\parallel
K_r\parallel_1$ satisfies 
$$\parallel K_r\parallel_1=O(r^{-1})$$
as $r\to\infty$. Thus
\begin{equation}\label{7.24}
K_{r,1}:=(\Id-\Pi_1) K_r(\Id-\Pi_1)
\end{equation}
is also a trace class operator with trace norm satisfying
\begin{equation}\label{7.25}
\parallel K_{r,1}\parallel_1=O(r^{-1}),\quad r\to\infty.
\end{equation}
Furthermore, $B(r)$ and $C(r)$ are finite rank operators with
\begin{equation}\label{7.26}
\parallel B(r)\parallel_1,\;\parallel C(r)\parallel_1=O(r^{-1}).
\end{equation}
Finally, $A(r)$ is a linear operator in the finite dimensional vector space
$L_1$ whose norm is also $O(r^{-1})$. This operator can be described more
explicitely as follows. First note that $L_1\subset \ker\Delta_Y\oplus
\ker\Delta_Y$ and hence we can replace $\Pi_1$ by the orthogonal projection
$\Pi_2$ of $\ker\Delta_Y\oplus\ker\Delta_Y$ onto $L_1$. Let 
$(V_1^+\cap V_2^+)^\perp_i\subset V_i^+$ denote the orthogonal complement
of $V_1^+\cap V_2^+$ in $V_i^+$, $i=1,2$, and let 
$$P_i:\ker\Delta_Y\to (V_1^+\cap V_2^+)_i^\perp$$ be the
orthogonal projection of $\ker\Delta_Y$ onto $(V_1^+\cap V_2^+)^\perp_i$.
Then $\Pi_2=(P_1,P_2)$ and by (\ref{7.22}) it follows that
\begin{equation*}
A(r)=\frac{1}{2r}\begin{pmatrix}P_1 & 0\\0 & P_2\end{pmatrix}\circ
\begin{pmatrix}\Id & -\Id\\ -\Id& \Id\end{pmatrix}\circ
\begin{pmatrix}P_1 & 0\\0 & P_2\end{pmatrix}=\frac{1}{2r}
\begin{pmatrix}P_1 & -P_1P_2\\-P_2P_1 & P_2\end{pmatrix}.
\end{equation*}
Regarded as operator in $(V_1^+\cap V_2^+)^\perp_1\oplus
(V_1^+\cap V_2^+)^\perp_2$, we get 
\begin{equation}\label{7.27}
A(r)=\frac{1}{2r}\begin{pmatrix}\Id & -P_1\\ -P_2 &\Id\end{pmatrix}.
\end{equation}
Suppose that $(\phi,\psi)\in L_1$ is in the kernel of $A(r)$. Then it follows
that 
$$\phi=P_1\psi,\quad \psi=P_2\phi.$$
Since $\phi\in V_1^+$ and $\psi\in V_2^+$, it follows that $\phi,\psi\in
V_1^+\cap V_2^+$ and therefore $\phi=\psi=0$.
Thus $A(r)$ is invertible and its norm satisfies
\begin{equation}\label{7.28}
\parallel A(r)\parallel =c r^{-1},\quad r>0.
\end{equation}
for some constant $c>0$. Let
$$S\colon (V_1^+\cap V_2^+)_1^\perp\oplus (V_1^+\cap V_2^+)_2^\perp\to
(V_1^+\cap V_2^+)_1^\perp\oplus (V_1^+\cap V_2^+)_2^\perp
$$
denote the restriction of the operator
\begin{equation*}
\begin{pmatrix}\Id & -P_1\\-P_2 & \Id\end{pmatrix}\colon \ker\Delta_Y\oplus
\ker\Delta_Y\to\ker\Delta_Y\oplus\ker\Delta_Y
\end{equation*}
to the subspace $(V_1^+\cap V_2^+)_1^\perp\oplus (V_1^+\cap V_2^+)_2^\perp$.
Set 
\begin{equation}\label{7.28b}
h:=\dim V_1^+ +\dim V_2^+ -2\dim V_1^+\cap V_2^+ \quad\mathrm{and}\quad 
 h_{12}:=\dim V_1^+\cap V_2^+. 
\end{equation}
\begin{lem}\label{l7.9}
We have 
$$\lim_{r\to\infty}r^{h+h_{12}}\det R_r=2^{-h}\det(S)\det R_{1,\infty} 
\det R_{2,\infty}.$$
\end{lem}
\begin{proof}
Let
\begin{equation*}
T_0(r)=\begin{pmatrix}A(r) & 0\\0 & D(r)\end{pmatrix}.
\end{equation*}
Since $A(r)$ is an invertible operator in a finite-dimensional vector space
and $D(r)$ is invertible for $r\ge r_0$, it follows that $T_0(r)$ is
invertible for $r\ge r_0$ and
$$\det T_0(r)=\det A(r)\det D(r).$$
Let
\begin{equation*}
T_1(r)=\begin{pmatrix}0 & B(r)\\ C(r)& 0\end{pmatrix}.
\end{equation*}
Then $T_1(r)$ is a trace class operator with 
$\parallel T_1(r)\parallel_1=O(r^{-1})$ as $r\to\infty$, and
$$R_r=T_0(r)+T_1(r)$$
Moreover $T_0(r)+tT_1(r)$ is invertible for $0\le t\le 1$ and $r\ge r_0$. Put
\begin{equation*}
T_2(r)=T_1(r)T_0(r)^{-1}.
\end{equation*}
Then for $r\ge r_0$ we get
\begin{equation}\label{7.28a}
\begin{split}
\log\det R_r-\log\det T_0(r)&=\int_0^1\frac{d}{dt}\log\det(T_0(r)+tT_1(r))
\;dt\\
&=\int_0^1\Tr\left(T_1(r)(T_0(r)+tT_1(r))^{-1}\right)\;dt\\
&=\int_0^1\Tr\left(T_2(r)(\Id+tT_2(r)\right)^{-1}\;dt.
\end{split}
\end{equation}
Set 
$$\widetilde B(r)=B(r)D(r)^{-1},\quad \widetilde C(r)=C(r)A(r)^{-1}.$$
Using the definition of $T_2(r)$, we get
\begin{equation*}
T_2(r)=\begin{pmatrix}0 & \widetilde B(r)\\ 
\widetilde C(r)& 0
\end{pmatrix}
\end{equation*}
By (\ref{7.26}) and (\ref{7.28}) it follows that
\begin{equation}\label{7.29}
\parallel \widetilde B(r)\parallel=O(r^{-1}),\quad \parallel\widetilde C(r)
\parallel =O(1)
\end{equation}
as $r\to\infty$. Thus
\begin{equation*}
T_2(r)^2=\begin{pmatrix}\widetilde B(r)\widetilde C(r) & 0\\
{} & {}\\
0 & \widetilde C(r)\widetilde B(r)\end{pmatrix}
\end{equation*}
and by (\ref{7.29}) we have
\begin{equation*}
\parallel T_2(r)^2\parallel=O(r^{-1}),\quad r\to\infty.
\end{equation*}
Let $r_1>0$ be such that
$$\parallel T_2(r)^2\parallel\le \frac{1}{2}$$
for $r\ge r_1$. Then
\begin{equation*}
\sum_{k=0}^\infty T_2(r)^k=(\Id+T_2(r))\sum_{k=0}^\infty T_2(r)^{2k}
\end{equation*}
is absolutely convergent and hence, $\Id+tT_2(r)$ is invertible for
$0\le t\le 1$ and $r\ge r_1$ with
\begin{equation*}
\left(\Id+tT_2(r)\right)^{-1}=\sum_{k=0}^\infty t^kT_2(r)^k.
\end{equation*}
Moreover it follows that
\begin{equation*}
\left(\Id+tT_2(r)\right)^{-1}=\Id+tT_2(r)+O(r^{-1}),\quad r\to\infty.
\end{equation*}
Thus
$$T_2(r)(\Id+tT_2(r))^{-1}=T_2(r)+tT_2(r)^2+O(r^{-1})=T_2(r)+O(r^{-1}).$$
Since $\Tr(T_2(r))=0$, we get
$$\Tr\bigl(T_2(r)(\Id+tT_2(r))^{-1}\bigr)=O(r^{-1}).$$
Together with (\ref{7.28a}) this implies
$$\big|\log\det R_r -\log\det T_0(r)\big|\le C r^{-1}.$$
Hence we get
$$\frac{\det (R_r|L^\perp)}{\det T_0(r)}=1+O(r^{-1}),\quad r\to\infty.$$
As observed above, $\det T_0(r)=\det A(r)\det D(r)$. Using the definition
of $D(r)$ by (\ref{7.23a}) and that $R_\infty|(\ker R_\infty)^\perp$ is
invertible, it follows as in Lemma 6.5 that
$$\lim_{r\to\infty} \det D(r)=\det R_{1,\infty}\det R_{2,\infty}.$$
Let $h$ and $h_{12}$ be defined by (\ref{7.28b}). Note that 
$$h=\dim(V_1^+\cap V_2^+)_1^\perp +\dim(V_1^+\cap V_2^+)_2^\perp.$$
Then by definition of $A(r)$ 
$$\det A(r)=(2r)^{-h}\det S.$$
So combined with (\ref{7.22a})  we get
$$\lim_{r\to\infty}r^{h+h_{12}}\det R_r=2^{-h}\det (S)\det R_{1,\infty}
\det R_{2,\infty}.$$
\end{proof}

Next we express $\det(S)$ in terms of the scattering matrices $S_1(0)$ and 
$S_2(0)$. Let $V=\ker \Delta_Y$ and set 
$$V_2 =V\ominus((V_1^+\cap V_2^+)\oplus (V_1^-\cap V_2^-)).$$ 

\begin{lem}\label{l7.10}
Let $C_{12}=S_1(0)S_2(0)|V_2$. We have 
$$\det(S)=\det\left((\Id-C_{12})/2 \right).$$
\end{lem}
\begin{proof}
First we consider the following  special case: Assume that 

\begin{enumerate}
\item $V_1^+\cap V_2^+=\{0\},\quad V_1^-\cap V_2^-=\{0\},$
\item $\dim V=2p$ and $ \dim V_i^+=\dim V_i^-=p,\;i=1,2,$
\item $P_1^+\colon V_2^+\to V_1^+$ is an isomorphism.
\end{enumerate}

Let $e_1,...,e_{2p}$ be an orthonormal basis of
$\ker\Delta_Y$ such that $e_1,...,e_p$ is an orthonormal basis of $V_1^+$ and
$e_{p+1},...,e_{2p}$ is an orthonormal basis of $V_1^-$. Let $f_1,...,f_p\in
V_2^+$  be such that
$$P_1^+(f_i)=e_i,\quad i=1,...,p.$$
Then there exists a symmetric matrix $A=(a_{ij})\in\GL(p,\R)$ such
 that
\begin{equation*}
f_i=e_i+\sum_{j=1}^p a_{ij}e_{p+j},\quad i=1,...,p.
\end{equation*}
Let $A^{-1}=(b_{kl})$ and put
\begin{equation*}
f_{p+k}=e_k+\sum_{l=1}^p b_{kl}e_{p+k},\quad k=1,...,p.
\end{equation*}
Then $\langle f_i,f_{p+j}\rangle =0$, $i,j=1,...,p$. Thus $f_{p+j}\in V_2^-$,
$j=1,...,p$. Furthermore $P_1^+(f_{p+j})=e_j$.  Thus $f_{p+1},...,f_{2p}$ is a
basis of $V_2^-$. By definition the  matrix $T$ which transforms the basis 
$(e_1,...,e_{2p})$ into $(f_1,...,f_{2p})$ is given by
\begin{equation*}
T=\begin{pmatrix}\Id & A\\ \Id &-A^{-1}\end{pmatrix}.
\end{equation*}
Since $A$ is symmetric, it follows that $(A^2+\Id)$ is invertible and one
immediately verifies that the inverse of $T$ is given by
\begin{equation*}
T^{-1}=\begin{pmatrix}(A^2+\Id)^{-1} & A^2(A^2+\Id)^{-1}\\ {} & {}\\
A(A^2+\Id)^{-1} &-A(A^2+\Id)^{-1}\end{pmatrix}.
\end{equation*}
Now note that the matrix of $S_1(0)S_2(0)$ with respect to the basis
$(e_1,...,e_{2p})$ is given by
\begin{equation}\label{7.30}
\begin{pmatrix}\Id & 0\\ 0 & -\Id\end{pmatrix}\circ 
 T^{-1}\circ \begin{pmatrix}\Id & 0\\ 0 & -\Id\end{pmatrix}\circ T
\end{equation}
Hence the matrix of $\Id-S_1(0)S_2(0)$ in the basis
$(e_1,...,e_{2p})$ is equal to
\begin{equation*}
\begin{pmatrix} 2A^2(A^2+\Id)^{-1} & -2A(A^2+\Id)^{-1}\\{} &{}\\
2A(A^2+\Id)^{-1} & 2A^2(A^2+\Id)^{-1}\end{pmatrix}.
\end{equation*}
This implies 
\begin{equation*}
\det(\Id-S_1(0)S_2(0))=2^{h_Y}\det(A^2)\det(A^2+\Id)^{-1}.
\end{equation*}

On the other hand $P_2^+=1/2(\Id+S_2(0))$. So it follows from (\ref{7.30})
that in the basis $(e_1,...,e_{2p})$, $P_2^+$ is represented by the matrix
\begin{equation*}
\begin{pmatrix}
(A^2+\Id)^{-1} & A(A^2+\Id)^{-1}\\{} &{}\\
A(A^2+\Id)^{-1} & A^2(A^2+\Id)^{-1}\end{pmatrix}.
\end{equation*}
Thus, with respect to the bases $(e_1,...,e_p)$ and $(f_1,...,f_p)$, the 
operator $P_2^+:V_1^+\to V_2^+$ is represented by the matrix $(A^2+\Id)^{-1}$.
Hence the matrix of $S$ with respect to the basis $(e_1,...,e_p,f_1,...,f_p)$
is given by
\begin{equation*}
\begin{pmatrix}
\Id & -\Id\\
-(A^2+\Id)^{-1} & \Id\end{pmatrix}.
\end{equation*}
Thus
\begin{equation*}
\det(S)=\det\begin{pmatrix}
\Id-(A^2+\Id)^{-1} & 0\\
-(A^2+\Id)^{-1} & \Id\end{pmatrix}
=\det(A^2)\det(A^2+\Id)^{-1}.
\end{equation*}

Next we reduce the general case to this special one. If we restrict $S_1(0)$
 and
$S_2(0)$ to $V_2$, it follows immediately that we can assume condition 1). Now
suppose that
$$\dim V_2^+\le \dim V_1^+\quad\text{and}\quad P_1^+:V_2^+\to V_1^+$$
is injective. Let $W_1:=P_1^+(V_1^+)$ and let  $W_2\subset V_1^+$ denote the
orthogonal complement of $W_1$ in $V_1^+$. We claim that $W_2\subset V_2^-$. 
To prove this claim let $w\in W_2$ and $v\in V_2^+$ be given. Write
$v=v_1+v_2, \quad v_1\in V_1^+,\; v_2\in V_1^-$. By definition we have
$\langle w,v_1\rangle=0$. Since $w\in W_2\subset V_1^+$, we have
$\langle w,v_2\rangle=0$. Thus $\langle w,v\rangle=0$, which shows that
$W_2$ is orthogonal to $V_2^+$, and hence $W_2\subset V_2^-$. Now
\begin{equation}\label{7.33a}
S_1(0)|W_2=\Id,\quad S_2(0)|W_2=-\Id.
\end{equation}
Thus $S_1(0)S_2(0)|W_2=-\Id$. Let
$$\tilde V=V\ominus W_2.$$
Then by (\ref{7.33a}), $\tilde V$ is an invariant subspace for $S_1(0)$ and
$S_2(0)$. Let $\tilde S_i=S_i(0)|\tilde V$, $i=1,2$. Then
\begin{equation*}
\Id-S_1(0)S_2(0)=\begin{pmatrix} 2\Id & 0\\ 0& \Id-\tilde S_1\tilde S_2
\end{pmatrix}.
\end{equation*}
Hence we get
\begin{equation}\label{7.33b}
\det(\Id-S_1(0)S_2(0))=2^{\dim W_2}\det(\Id-\tilde S_1\tilde S_2).
\end{equation}
Let $\tilde V_i^\pm\subset \tilde V$ be the $\pm1$-eigenspaces of $\tilde
S_i$, $i=1,2$. Then it follows that
$$\tilde V_1^+=P_1^+(V_2^+)=W_1,\quad \tilde V_2^+=V_2^+.$$
In particular, $\tilde P_1^+:\tilde V_2^+\to \tilde V_1^+$ is an isomorphism.
Thus $\dim \tilde V_1^\pm=\tilde V_2^\pm$. Since $\tilde V_1^\pm\cap \tilde
V_2^\pm=\{0\}$, it follows that $\dim\tilde V_i^+=1/2\dim\tilde V$. Thus
conditions 2) and 3)  are also satisfied and hence, by the first part of the
proof we get
\begin{equation}\label{7.33c}
\det(\Id-\tilde S_1\tilde S_2)
=2^{\dim\tilde V}\det\begin{pmatrix} \Id & -\tilde P_1^+\\ -\tilde P_2^+&
    \Id \end{pmatrix}.
\end{equation}
Finally note that with respect to the decomposition $V_1^+=W_1\oplus W_2$,
$$P_1^+\colon V_2^+\to V_1^+\quad \text{and}\quad P_2^+\colon V_1^+\to V_2^+$$
are of the form
$$P_1^+=(\tilde P_1^+,0),\quad P_2^+=\tilde P_2^+\oplus 0.$$
Hence
\begin{equation*}
\begin{pmatrix}\Id & -P_1^+\\ -P_2^+ & \Id\end{pmatrix}=\begin{pmatrix}\Id &0& 
-\tilde P_1^+\\ 0 &\Id & 0\\ -\tilde P_2^+ & 0 &\Id\end{pmatrix},
\end{equation*}
which shows that
$$\det\begin{pmatrix}\Id_{V_1^+} & -P_1^+\\ -P_2^+ & \Id_{V_2^+}
\end{pmatrix}=\det\begin{pmatrix}\Id_{\tilde V_1^+} & -\tilde P_1^+\\ 
-\tilde P_2^+ & \Id_{\tilde V_2^+}\end{pmatrix}.$$
Together with (\ref{7.33b}) and (\ref{7.33c}), the lemma follows.
\end{proof}

Combining Proposition \ref{p6.4} with Corollary \ref{c7.3} and Lemmas
\ref{l7.7}, \ref{l7.9} and \ref{l7.10}, we obtain
\begin{equation}\label{7.37}
\begin{split}
\lim_{r\to\infty}r^{h}\frac{\det \Delta_{M_r}}{\det\Delta_{N_r,D}}&=
\prod_{i=1}^2\frac{\det(\Delta_{i,\infty},\Delta_0)}{\det
R_{i,\infty}}\lim_{r\to\infty} \frac{r^{h+h_{12}}
\det R_r}{r^{h_{12}}\det B_r}\\
&=2^{-h}\det\left((\Id-C_{12})/2\right)
\prod_{i=1}^2\det(\Delta_{i,\infty},\Delta_0). 
\end{split}
\end{equation}
Using Proposition (\ref{p5.1}), we get 
\begin{equation}\label{7.38}
\begin{split}
\det\Delta_{M_r}\sim r^{h_Y-h}&e^{-r\xi^\prime_Y(0)}2^{2h_Y-h}
(\det\Delta_Y)^{-1/2}\\
&\cdot\det\left((\Id-C_{12})/2\right)
\prod_{i=1}^2\det(\Delta_{i,\infty},\Delta_0).
\end{split}
\end{equation}
As an example, we consider the case of a closed surface $M.$ Let
$$M_L=M_1\cup_Y([0,L]\times Y)\cup_YM_2,\quad Y=\R/\Z,\;L>0.$$
Then
$$\zeta_Y(s)=(2\pi)^{-2s}2\zeta(2s),$$
where $\zeta(s)$ denotes the Riemann zeta function.
Recall that
$$\zeta(-1)=-\frac{1}{12},\quad \zeta(0)=- \frac{1}{2},\quad
\zeta^\prime(0)=-\frac{1}{2}\log 2\pi.$$
Since $^{\cdot}\Gamma(s-1/2)$ and $\zeta(2s)$ are analytic at $s=0$, we get
$$\xi^\prime_Y(0)=\frac{1}{\sqrt{\pi}}\frac{d}{ds}\left(\frac{\Gamma
(s-1/2)}{\Gamma(s)}\zeta_Y(s-1/2)\right)\Big|_{s=0}
=-\frac{\sqrt{\pi}}{3}\frac{d}{ds}
\left(\frac{\Gamma(s-1/2)}{\Gamma(s)}\right)\Big|_{s=0}=\frac{2}{3}\pi.$$
Similarly
$$\zeta^\prime_Y(0)=-4\log(2\pi)\zeta(0)+4\zeta^\prime(0)=0.$$
Thus
$$\det\Delta_Y=e^{-\zeta{^\prime_Y}(0)}=1.$$
Furthermore  note that $h_Y=h_{12}=1$, $h=0$, and 
$\det\left((\Id-C_{12})/2\right)=1/2$.  
 Inserting this into (\ref{7.38}), we get
\begin{equation}\label{7.39}
\det\Delta_{M_L}\sim 2L e^{-\pi L/3}\det(\Delta_{1,\infty},\Delta_0)
\cdot\det(\Delta_{2,\infty},\Delta_0)
\end{equation}
as $L\to\infty.$ Bismut and Bost proved in \cite{BB} that
$\det\Delta_{M_L}\sim c L e^{-\pi L/3},$ $L\to\infty,$ with some constant
$c.$ Our result expresses the constant $c$ explicitely as 
$c=2\det(\Delta_{1,\infty},\Delta_0)\cdot\det(\Delta_{2,\infty},\Delta_0)$.

Next consider a compact Riemannian manifold $(X_0,g)$ with boundary $Y$ as at
the end of the previous section. We assume that the connection $\nabla^E$ is a
product on the collar neighborhood $N=(-\epsilon,0]\times Y$ of $Y$ in $X_0$. 
By (\ref{6.17}) and Corollary \ref{c7.3} we have
\begin{equation}\label{7.40}
%\begin{split}
\log\det \Delta_{X_r,D}=\log\det\Delta_{N_r,D}+\log\det R_r 
-\log\det R_\infty +\det(\Delta_\infty,\Delta_0).
%\end{split}
\end{equation}
Furthermore by Lemma \ref{l7.2} we have $\ker R_\infty=V^+$. By (\ref{6.17a}) 
it follows that $\ker R_\infty$ is invariant under $L_r$ and hence under
$R_r$, and
$$R_r|_{\ker R_\infty}=\frac{1}{r}\Id.$$
Let $h^+=\dim V^+$. Then 
$$\det R_r=r^{-h^+}\det\left(R_r|(\ker R_\infty)^\perp\right)$$
and by Lemma 4.1 of \cite{Le3} it follows that
$$\lim_{r\to\infty} r^{h^+}\det R_r=\det R_\infty.$$
Using (\ref{7.40}) we obtain
$$\lim_{r\to\infty}\frac{r^{h^+}\det\Delta_{X_r,D}}{\det\Delta_{N_r,D}}=
\det(\Delta_\infty,\Delta_0).$$
Together with Proposition \ref{p5.1} we get
$$\det\Delta_{X_r,D}\sim r^{-h^+ +h_Y}e^{-r\xi_Y^\prime(0)/2}2^{h_Y}
\left(\det\Delta_Y\right)^{-1/2}\det(\Delta_\infty,\Delta_0)$$
as $r\to\infty$. If we apply this to $\det\Delta_{M_{i,r},D}$ and use
(\ref{1.9}), we get
$$\lim_{r\to\infty}\frac{r^{h_Y-2h_{12}}\det\Delta_{M_r}}
{\det\Delta_{M_{1,r},D}\det\Delta_{M_{2,r},D}}=2^{-h}
\left(\det\Delta_Y\right)^{1/2}\det\left((\Id-C_{12})/2\right),$$
which proves Theorem \ref{th1.7}.


\begin{thebibliography}{HMM}
\bibitem[Ba]{Ba} C. B\"ar, {\it  Zero sets of solutions to semilinear elliptic
    systems of first order,}  Invent. Math.  {\bf 138} (1999), 183--202.
\bibitem[BB]{BB} J.-M. Bismut, J.-B. Bost, {\it Fibr\'es d\'eterminants,
m\'etriques de Quillen et d\'eg\'en\'erescence de courbes}, Acta Math. {\bf
165} (1990), 1-103. 
\bibitem[BFK]{BFK} D. Burghelea, L. Friedlander, and T. Kappeler, 
{\it Meyer-Vietoris type formula for determinants of elliptic differential
operators}, J. Funct. Analysis {\bf 107} (1992), 34--65.
\bibitem[Ca]{Ca} G. Carron, {\it D\'eterminant relatif et fonction Xi}, 
 Amer. J. Math.  {\bf 124} (2002), 307--352.
\bibitem[Do]{Do} H. Donnelly, {\it Eigenvalue estimates for certain noncompact
manifolds}, Michigan Math. J. {\bf 31} (1984), 349--357.
o\bibitem[Gi]{Gi} P. Gilkey, {\it The spectral geometry of a riemannian
manifold}, J. Differential Geometry {\bf 10} (1975), 601--618.
\bibitem[Gu]{Gu} L. Guillop\'e, {\it Th\'eorie spectrale de quelques 
vari\'et\'es \`a bouts}, Ann.scient. \'Ecole Norm. Sup., $4^e$ s\'erie, t.
{\bf 22} (1989), 137-160.
\bibitem[Ha]{Ha} A. Hassell,  {\it Analytic surgery and analytic torsion.}
  Comm. Anal. Geom.  {\bf 6} (1998), 255--289. 
\bibitem[HMM]{HMM} A. Hassell, R. Mazzeo, R.B. Melrose, {\it analytic surgery
and the accumulation of eigenvalues}, Comm. Analysis and Geometry. {\bf 3}
(1995), 115-222. 
\bibitem[HZ]{HZ} A. Hassell, S. Zelditch, {\it Determinants of Laplacians in
exterior domains}, IMRN {\bf 18} (1999), 971--1004.
\bibitem[La]{La} S. Lang, {\it Elliptic functions}, Addison-Wesley, 1973.
\bibitem[Le1]{Le1} Y. Lee, {\it Mayer-Vietoris formula for the determinant of
a Laplace operator on an even-dimensional manifold}, Proc. Amer. Math. Soc. 
{\bf 123} (1995), 1933--1940.
\bibitem[Le2]{Le2} Y. Lee, {\it Mayer-Vietoris formula for determinants of
elliptic operators of Laplace-Beltrami type}, Diff. Geometry and its 
Applications {\bf 7} (1997), 325--340.
\bibitem[Le3]{Le3} Y. Lee, {\it Burghelea-Friedlander-Kappeler's gluing
formula for the zeta-determiannt and its applications to the adiabatic
decompositions of the zeta-determinant and the analytic torsion}, Preprint
2003, arXiv:math.DG/0304250.
\bibitem[Le4]{Le4} Y. Lee, {\it Asymptotic expansion of the zeta-determinant of
an invertible Laplacian operator on a streched manifold}, to appear in
 Contemporary Math.
\bibitem[LP]{LP} P. Loya, J. Park, {\it Decomposition of the 
$\zeta$-determinant for the Laplacian on manifolds with cylindrical end},
Preprint, 2003.
\bibitem[MM] {MM} R. Mazzeo and R.B. Melrose, {\it Analytic surgery and the eta
  invariant}.  Geom. Funct. Anal.  {\bf 5}(1995), 14--75.  
\bibitem[Mu1]{Mu1} W. M\"uller, {\it Relative zeta functions, relative
determinants and scattering theory}, Commun. Math. Phys. {\bf 192} (1998),
309--347.
\bibitem[Mu2]{Mu2} W. M\"uller, {\it $L^2$-index and resonances}, In: 
''Geometry and Analysis'', LNM 1339, Springer-Verlag, 1988, pp. 203-211.
\bibitem[Mu3]{Mu3} W. M\"uller, {\it Manifolds with cusps of rank one},
Lecture Notes in Math. {\bf 1244}, Springer-Verlag, 1987.
\bibitem[Mu4]{Mu4} W. M\"uller, {\it Eta invariants and manifolds with
boundary}, J. Differential Geom. {\bf 40} (1994), 311-377.
\bibitem[PW1]{PW1} J. Park, K.P. Wojciechowski, 
{\it Adiabatic decomposition of the $\zeta$-determinant of the Dirac Laplacian
  I. The case of invertible tangential operator}. With an appendix by Y. Lee. 
 Comm. Partial Differential
Equations {\bf 27} (2002), 1407-1435.
\bibitem[PW2]{PW2} J. Park, K.P. Wojciechowski, {\it Scattering theory and
    adiabatic decomposition of the $\zeta$-determinant of the Dirac
    Laplacian.}  Math. Res. Lett. {\bf 9} (2002), 17--25.
\bibitem[PW3]{PW3} J. Park, K.P. Wojciechowski, {\it Adiabatic decomposition
    of the $\zeta$-determinant and Dirichlet to Neumann operator}, Preprint
2003, arXiv:math.DG/0301170.  
\end{thebibliography}
\end{document}